\documentclass[12pt]{article}
\usepackage[utf8]{inputenc}
\usepackage[english]{babel}
\usepackage{pgf,tikz,tikz-cd}

\definecolor{LightGray}{rgb}{0.85,0.85,0.85}
\usepackage{amsmath,mathtools,amsfonts,amscd,amssymb}
\usepackage[all]{xy}
\usepackage{makeidx}
\usepackage{color}
\usepackage{graphicx}
\usepackage[active]{srcltx} 
\newtheorem{thm}{Theorem}[section]
\newtheorem{cor}[thm]{Corollary}
\newtheorem{lem}[thm]{Lemma}
\newtheorem{prop}[thm]{Proposition}
\newtheorem{deft}[thm]{Definition}
\newtheorem{rek}[thm]{Remark}
\newtheorem{conj}{Conjecture}
\numberwithin{equation}{section}
\let\noi=\noindent

\let\sur=\overline

\def\N{\mathbb{N}} 
\def\Z{\mathbb{Z}} 
\def\R{\mathbb{R}} 
\def\C{\mathbb{C}} 
\def\Q{\mathbb{Q}} 
\def\F{\mathbb{F}}
\def\sper{\mbox{Sper}}

\def\gr{\mbox{gr}}
\def\init{\mbox{in}}
\def\notin{\mbox{$\in$ \hspace{-.8em}/}} 
\def\sgn{\mbox{sgn}} 

\makeindex

\title{ON THE STRONG SEPARATION CONJECTURE}
\author{F. LUCAS $^\dagger$ $^a$, D. SCHAUB $^a$,  M. SPIVAKOVSKY $^b$ \\  \\ $^a$ Larema, Université Angers, CNRS UMR 6093, 2, bd Lavoisier, \\ 49045 Angers Cedex, France \\ $^b$ 
	Institut de Math\'ematiques de Toulouse, CNRS UMR 5219, UNAM\\
	Universit\'e Paul Sabatier, 118, route de Narbonne \\
	31062 Toulouse Cedex 9, France.}

\date{}

\begin{document}
\maketitle

\footnotetext[2]{ \quad F. Lucas passed away in April 2016}

\abstract{This paper contains a partial result on the Pierce--Birkhoff conjecture on piecewise polynomial functions defined by a finite collection $\left\{f_1,\dots,f_r\right\}$ of polynomials. In the nineteen eighties, generalizing the problem from the polynomial ring to an artibtrary ring $\Sigma$, J. Madden proved that the Pierce--Birkhoff conjecture for $\Sigma$ is equivalent to a statement about an arbitrary pair of points $\alpha,\beta\in\sper\ \Sigma$ and their separating ideal $<\alpha,\beta>$; we refer to this statement as the \textbf{local Pierce-Birkhoff conjecture} at $\alpha,\beta$. In \cite{LMSS} we introduced a slightly stronger conjecture, also stated for a pair of points $\alpha,\beta\in\sper\ \Sigma$ and the separating ideal $<\alpha,\beta>$, called the \textbf{Connectedness conjecture}, about a finite collection of elements $\left\{f_1,\ldots,f_r\right\}\subset\Sigma$. In the paper \cite{LMSS3} we introduced a new conjecture, called the \textbf{Strong Connectedness conjecture}, and proved that the Strong Connectedness conjecture in dimension $n-1$ implies the Strong Connectedness conjecture in dimension $n$ in the case when $ht(<\alpha,\beta>)\le n-1$. 

The Pierce-Birkhoff Conjecture for $r=2$ is equivalent to the Connectedness Conjecture for $r=1$; this conjecture is called the Separation Conjecture. The Strong Connectedness Conjecture for $r=1$ is called the Strong Separation Conjecture. In the present paper, we fix a polynomial $f \in R[x,z]$ where $R$ is a real closed field and $x=(x_1,\ldots,x_n),z$ are $n+1$ independent variables. We define the notion of two points $\alpha,\beta \in \sper\ R[x,z]$ being in \textbf{good position} with respect to $f$. The main result of this paper is a proof of the Strong Separation Conjecture in the case when $\alpha$ and $\beta$ are in good position with respect to
$f$.}
\bigskip

\noi\textit{Dedicated to Professor Felipe Cano on the occasion of his sixtieth birthday.}

\section{Introduction}

All the rings in this paper will be commutative with 1. Let $R$ be a real closed field. Let $x=(x_1,\ldots,x_n)$ and let $z$ be a single variable. Let $B=R[x]$. We will use the notation $A:=R[x,z]=B[z]$.

Throughout the paper, by \textbf{connectedness} we mean \textbf{semi-algebraic connectedness} (sometimes called definable connectedness, cf. \cite{BCR}, D\'efinition 2.4.5). If $R=\R$, a semi-algebraic subset of $R^n$ is connected if and only if it is semi-algebraically connected (\cite{BCR}, Theorem 2.4.5). 
\medskip 


The Pierce--Birkhoff conjecture asserts that any piecewise-polynomial function
\[
g:R^n\rightarrow R
\]
can be expressed as a maximum of minima of a finite family of polynomials in $n$ variables.

We start by giving a precise statement of the conjecture as it was
first stated by M. Henriksen and J. Isbell in the early nineteen
sixties (\cite{BiPi} and \cite{HenIsb}).
\begin{deft}\label{pw} A function $g:R^n\to R$ is said to be \textbf{piecewise
  polynomial} if $R^n$ can be covered by a finite collection of closed
  semi-algebraic sets $P_i$, $i\in\{1,\dots,s\}$ such that for each $i$ there exists a
  polynomial $g_i\in B$ satisfying
  $\left.g\right|_{P_i}=\left.g_i\right|_{P_i}$.
\end{deft}
Clearly, any piecewise polynomial function is continuous. Piecewise
polynomial functions form a ring, containing $B$, which is denoted by
$PW(B)$.\medskip

On the other hand, one can consider the (lattice-ordered) ring of all
the functions obtained from $B$ by iterating the operations of $\sup$
and $\inf$. Since applying the operations of sup and inf to
polynomials produces functions which are piecewise polynomial, this
ring is contained in $PW(B)$ (the latter ring is closed under $\sup$
and $\inf$). It is natural to ask whether the two rings coincide. The
precise statement of the conjecture is:
\begin{conj}\textnormal{\textbf{(Pierce-Birkhoff)}}\label{PB} If
  $g:R^n\to R$ is in $PW(B)$, then there exists a finite family of
  polynomials $g_{ij}\in B$ such that
  $f=\sup\limits_i\inf\limits_j(g_{ij})$ (in other words, for all
  $x\in R^n$, $f(x)=\sup\limits_i\inf\limits_j(g_{ij}(x))$).
\end{conj}
Here is a partial list of earlier papers devoted to the Pierce--Birkhoff conjecture: \cite{De}, \cite{LMSS}, \cite{LMSS2}, \cite{LMSS3}, \cite{Mad1}, \cite{Mah}, \cite{Mah2} \cite{Mar2} and \cite{W}.

The starting point of this paper is the abstract formulation of the conjecture in terms of the real
spectrum of $B$ and separating ideals  proposed by J. Madden in
1989 \cite{Mad1}.

For more information about the real spectrum, see \cite{BCR}; there is
also a brief introduction to the real spectrum and its relevance to
the Pierce--Birkhoff conjecture in the Introduction to \cite{LMSS}.
\medskip

\noi\textbf{Terminology}: If $\Sigma$ is an integral domain, the phrase
``valuation of $\Sigma$'' will mean ``a valuation of the field of fractions
 of $\Sigma$, non-negative on $\Sigma$''. Also, we will sometimes commit the
 following abuse of notation. Given a ring $\Sigma$, a prime ideal
 $\mathfrak p\subset \Sigma$, a valuation $\nu$ of $\frac \Sigma{\mathfrak p}$
 and an element $x\in\Sigma$, we will write $\nu(x)$ instead of
 $\nu(x\mod\mathfrak p)$, with the usual convention that
 $\nu(0)=\infty$, which is taken to be greater than any element of the
 value group.
\medskip

Let $K$ be a field and $\nu$ a valuation of $K$ with value group $\Gamma$. 
\medskip

\noi \textbf{Notation.} For $\gamma \in \Gamma$, let
\[
P_\gamma =\left\{g \in K \setminus \{0\}\left|\  \nu(g) \ge \gamma\right.\right\} \cup \{0\}
\]
and
\[
P_{\gamma_ +} = \left\{g \in K \setminus \{0\} \left|\ \nu(g) > \gamma\right.\right\}\cup \{0\}.
\]

Let $\Sigma$ be a subring of $K$. We define the graded algebra associated to $\nu$ to be
\[
\gr_{\nu}(\Sigma) = \bigoplus_{\gamma \in \Gamma} \frac{P_\gamma \cap\Sigma}{P_{\gamma_+}\cap\Sigma}.
\]
For all $f \in K$, $f \in P_\gamma \setminus P_{ \gamma_+}$, we denote by 
$\init_\nu(f)$ the natural image of $f$ in
\[
\frac{P_\gamma}{P_{\gamma_+}} \subset \gr_\nu(K).
\]
\medskip

For a point $\alpha \in \sper\ \Sigma$ we denote by $\mathfrak{p}_\alpha$ the support of $\alpha$. We let
$\Sigma[\alpha]=\frac{\Sigma}{\mathfrak{p}_\alpha}$ and let $\Sigma(\alpha)$ be the field of fractions of
$\Sigma[\alpha]$. We let $\nu_\alpha$ denote the valuation of $\Sigma(\alpha)$ associated to
$\alpha$ (\cite{LMSS2}, p. 264), $\Gamma_\alpha$ its the value group, $R_{\nu_\alpha}$ the valuation
ring, $k_\alpha$ its residue field and $\gr_\alpha(\Sigma)$ the graded algebra
associated to the valuation $\nu_\alpha$. For $f \in\Sigma$ with
$\gamma=\nu_\alpha(f)$, let $\init_\alpha f$ denote the natural image of $f$ in
$\frac{P_\gamma}{P_{\gamma_+}}$. For an ordered field $k$, we will denote by $\sur{k}_r$ a real closure of $k$ and by $\sur k$ an algebraic closure of $\sur{k}_r$. Typically we will work with ordered fields of the form $k=\Sigma(\alpha)$. Usually we will write $\sur{k}_r$ (respectively, $\bar k$) to indicate that we have fixed a real closure $\sur{k}_r$ of $k$ (resectively, an algebraic closure $\bar k$ or $\sur{k}_r$) once and for all.
\medskip

Next, we recall our generalization of the notion of piecewise polynomial functions and the Pierce--Birkhoff conjecture from polynomials to arbitrary rings (\cite{LMSS2}, Definition 8). Let $\Sigma$ be a ring.
\begin{deft}
	Let
	\[
	g :\sper\ \Sigma\to \coprod_{\gamma \in \sper\ \Sigma} \Sigma(\gamma)
	\]
	be a map such that, for each $\gamma \in \sper\ \Sigma$ we have $g(\gamma) \in\Sigma(\gamma)$. We
	say that $g$ is \textbf{piecewise polynomial}
	(denoted by $f\in PW(\Sigma)$) if there exists a
	covering of $\sper\ \Sigma$ by a finite family $(S_i)_{i \in\{1,\dots,s\}}$ of
	constructible sets, closed in the spectral topology, and a family
	$(g_i)_{i \in\{1,\dots,s\}}$, $g_i \in\Sigma$ such that, for each $\gamma \in S_i$,
	$g(\gamma)=g_i(\gamma)$.
	
	We call $g_i$ a \textbf{local representative} of $g$ at $\gamma$ and denote it
	by $g_\gamma$ ($g_\gamma$ is not, in general, uniquely determined by
	$g$ and $\gamma$; this notation means that one such local
	representative has been chosen once and for all).
\end{deft}

Note that $PW(\Sigma)$ is naturally a lattice ring: it is equipped
with the operations of maximum and minimum. Each element of $\Sigma$
defines a piecewise polynomial function. In this way we obtain a
natural injection $\Sigma\subset PW(\Sigma)$.

\begin{deft}
	A ring $\Sigma$ is a \textbf{Pierce--Birkhoff ring} if, for each $g \in PW(\Sigma)$, there
	exists a finite collection $\{g_{ij}\} \subset\Sigma$ such that
	$g=\sup\limits_i\inf\limits_j g_{ij}$.
\end{deft}
\begin{conj}\textnormal{\textbf{(the Pierce--Birkhoff conjecture for regular rings)}} A regular ring $\Sigma$ is a  Pierce-Birkhoff ring.
\end{conj}

J.J. Madden reduced the Pierce--Birkhoff Conjecture to a purely local statement about separating ideals and the real spectrum. Namely, he introduced 

\begin{deft} \cite{Mad1} Let $\Sigma$ be a ring. For $\alpha,\beta\in\mbox{Sper}\ \Sigma$,
  the \textbf{separating ideal} of $\alpha$ and $\beta$, denoted by
  $<\alpha,\beta>$, is the ideal of $\Sigma$ generated by all the elements
  $f\in\Sigma$ that change sign between $\alpha$ and $\beta$, that is,
  all the $f$ such that $f(\alpha)\geq0$ and $f(\beta)\leq0$.
\end{deft}

\begin{deft} \label{def:localPB}
A ring $\Sigma$ is \textbf{locally Pierce-Birkhoff at} $\alpha, \beta$ if the
following condition holds. Let $g$ be a piecewise  polynomial
  function, let $g_\alpha \in\Sigma$ be a local representative of $g$ at
  $\alpha$ and $g_\beta \in\Sigma$ a local representative of $g$ at
  $\beta$. Then $g_\alpha-g_\beta \in <\alpha,\beta>$.
\end{deft}

The statement that, for a certain class $\mathcal{X}$ of rings, every ring $\Sigma\in \mathcal{X}$ is locally Pierce-Birkhoff for 
all $\alpha, \beta \in \sper\ \Sigma$ will be denoted by PB($\mathcal{X}$). If $\mathcal{X}$ consists of only one ring $\Sigma$, we will write PB($\Sigma$) instead of PB($\mathcal{X}$).

\begin{thm} (J. Madden \cite{Mad1}) The ring $\Sigma$ is Pierce-Birkhoff if and only if PB($\Sigma$) holds.
\end{thm}

In \cite{LMSS}, we introduced
\begin{deft}\textnormal{\textbf{(the Connectedness property)}}\label{conn} Let $\Sigma$ be a ring and
\[
\alpha,\beta\in\sper\ \Sigma.
\]
We say that $\Sigma$ has the \textbf{Connectedness property} at $\alpha$ and $\beta$ if for any finite collection $f_1,\dots,f_s$ of elements of $\Sigma \; \setminus<\alpha,\beta>$ there exists a connected set
\[  
C\subset\sper\ \Sigma
\]
such that $\alpha,\beta\in C$ and
  $C\cap\{f_i=0\}=\emptyset$ for $i\in\{1,\dots,s\}$ (in other words,
  $\alpha$ and $\beta$ belong to the same connected component of the
  set $\sper\ \Sigma\setminus \{f_1\dots f_s=0\}$).
\end{deft}
The statement that, for a certain class of rings $\mathcal{X}$, every ring $\Sigma  \in \mathcal{X}$  has the Connectedness 
Property for all $\alpha, \beta \in \sper\ \Sigma$ will be denoted by CP($\mathcal{X}$).
\begin{conj}\textnormal{\textbf{(the Connectedness conjecture for regular rings)}}\label{connConj} Let $\mathcal X$ denote the class of all the regular rings. Then $CP(\mathcal X)$ holds.
\end{conj}

In the paper \cite{LMSS}, we proved that CP($\mathcal{X}$) implies PB($\mathcal{X}$) where $\mathcal X$ is the class of all the polynomial rings over $R$. The proof given in  \cite{LMSS} applies verbatim to show that CP($\mathcal{X}$) implies PB($\mathcal{X}$) for any class $\mathcal{X}$ of rings whatsoever.
One advantage of CP($\Sigma$) is that it is a statement about elements of $\Sigma$ which makes no mention of piecewise polynomial 
functions; in particular, if $\Sigma$ is a polynomial ring, CP($\Sigma$) is a statement purely about polynomials.
\medskip

For a field $k$ and an ordered group $\Gamma$, we will denote by $k\left[\left[t^\Gamma\right]\right]$, the ring of generalized power series over $k$ with exponents in $\Gamma$, that is, the ring formed by all the expressions of the form
$\sum\limits_{\gamma\in W}c_\gamma t^\gamma$, where $c_\gamma\in k$ and $\gamma$ and $W$ is a well ordered subset of the semigroup $\Gamma_+$ of non-negative elements of $\Gamma$. The ring $k\left[\left[t^\Gamma\right]\right]$ is equipped with a natural
 $t$-adic valuation with value group $\Gamma$. An order on the field $k$ induces an order on $k\left[\left[t^\Gamma\right]\right]$ in a natural way. 
\begin{deft} Let $\Sigma$ be a ring and $k$ an ordered field. A $k$\textbf{-curvette} on
  $\sper\ \Sigma$ is a homomorphism of the form
\[
\alpha:\Sigma\to k\left[\left[t^\Gamma\right]\right],
\]
where $\Gamma$ is an ordered group. A $k$\textbf{-semi-curvette} is a
$k$-curvette $\alpha$ together with
  a choice of the sign data $\sgn\ x_1,\dots,\sgn\ x_r\in\{+,-\}$, where
  $x_1,...,x_r$ are elements of $\Sigma$ whose $t$-adic values induce an
  $\F_2$-basis of $\Gamma/2\Gamma$.
\end{deft}

In \cite{LMSS2} we explained how to associate to a point $\alpha$ of $\sper\ \Sigma$
a $\sur{( k_\alpha)}_r$-semi-curvette. Conversely, given an
ordered field $k$, a $k$-semi-curvette $\alpha$ determines a prime
ideal $\mathfrak{p}_\alpha$ (the ideal of all the elements of $\Sigma$
which vanish identically on $\alpha$) and a total ordering on
$\frac\Sigma{\mathfrak{p}_\alpha}$ induced by the ordering of the ring
$k\left[\left[t^\Gamma\right]\right]$ of formal power series. Hence a $k$-semi-curvette determines a point in $\sper\ \Sigma$. Below, we will 
often describe points in the real spectrum by
specifying the corresponding semi-curvettes.
\medskip

We will use the following notation throughout the paper. We denote $\sqrt{<\alpha,\beta>}$ by $\mathfrak{p}$; set
$\mu_\alpha:=\nu_\alpha(<\alpha,\beta>)$ and $\mu_\beta:=\nu_\beta(<\alpha,\beta>)$.
\medskip

\noi Let $\mathcal{X}$ be a class of rings. We use the following notation:

CP$_n(\mathcal{X})$ means that CP holds for all the rings $\Sigma \in \mathcal{X}$ such that $\dim \Sigma \le n$;

PB$_n(\mathcal{X})$ means that PB holds for all the rings $\Sigma \in \mathcal{X}$ such that $\dim \Sigma \le n$;

CP$_{\le n}(\mathcal{X)}$ means that the Connectedness Property holds for all the rings $\Sigma \in \mathcal{X}$ and
$\alpha,\beta\in\sper\ \Sigma$ such that $ht(\mathfrak{p}) \le n$;

PB$_{\le n}(\mathcal{X})$ means that the local Pierce-Birkhoff Conjecture holds for all the rings $\Sigma \in \mathcal{X}$ and $\alpha,\beta \in \sper\ \Sigma$ such that $ht(\mathfrak{p}) \le n$.
\medskip

In the situation of CP$_n$ for a certain regular ring $\Sigma$ of dimension $n$, assume that $ht(\mathfrak{p}) < n$. If one wants to 
proceed by induction on $n$, a natural idea is to try to reduce CP$_n(\Sigma$) to CP$_{n-1}(\Sigma_\mathfrak{p})$. 

The difficulty with this approach is that the CP$_{n-1}$ cannot be 
applied directly. Indeed, let $f_1,\dots,f_s$ be as in the
CP and let $\Delta_\alpha\subset\Gamma_\alpha$ denote the 
greatest isolated subgroup not containing $\nu_\alpha(\mathfrak p)$. \\ The 
hypothesis $f_i\notin<\alpha,\beta>$ does not imply that
$f_i\notin<\alpha,\beta>\Sigma_{\mathfrak p}$: it may happen that
$\nu_\alpha(f_i)$ $< \mu_\alpha$,
$\nu_\alpha(f_i)-\nu_\alpha(\mathfrak p)\in\Delta_\alpha$ and so
$f_i\in<\alpha,\beta>\Sigma_{\mathfrak p}$, as shown by the Example below.
\smallskip

\noindent\textbf{Example.} Let $\Gamma= \Z^2_{lex}$. Let $\alpha,\beta \in \sper\ R[x,y,z]$ be given by the semi-curvettes
\begin{align}
x(t)&=t^{(0,3)}\label{eq:xt3}\\
y(t)&=t^{(0,4)}+bt^{(1,0)}\\
z(t)&=t^{(0,5)}+ct^{(1,1)}\label{eq:zt5},
\end{align}
where $b\in\{b_\alpha,b_\beta\}\subset\R$, 
$c\in\{c_\alpha,c_\beta\}\subset\R$ and $t^{(0,1)}>0, t^{(1,0)}>0$. The
constants $b_\alpha \neq b_\beta,c_\alpha \neq c_\beta$ will be specified later. Let
$f_1=xz-y^2$, $f_2=x^3-yz$, $f_3=x^2y-z^2$; consider the ideal
$(f_1,f_2,f_3)$. The most general common specialization of $\alpha,\beta$ is 
given by the semi-curvette
\begin{align}
x(t)&=t^3\\
y(t)&=t^4\\
z(t)&=t^5,
\end{align}
$t>0$. The corresponding point of $\sper\ R[x,y,z]$ has support $(f_1,f_2,f_3)$, so 
\[
\mathfrak p=\sqrt{<\alpha,\beta>}=(f_1,f_2,f_3).
\]
Let  $(x_\alpha(t),y_\alpha(t),z_\alpha(t))$ and 
$(x_\beta(t),y_\beta(t),z_\beta(t))$ 
be the semi-curvettes defining $\alpha$ and $\beta$ as in 
(\ref{eq:xt3})--(\ref{eq:zt5}). Let us calculate 
$f_i(x_\alpha(t),y_\alpha(t),z_\alpha(t))$ and 
$f_i(x_\beta(t),y_\beta(t),z_\beta(t))$. In the notation of 
(\ref{eq:xt3})--(\ref{eq:zt5}) we have
\begin{align}
f_1(x(t),y(t),z(t))&=(c-2b)t^{(1,4)}+\tilde f_1\\
f_2(x(t),y(t),z(t))&=-(c+b)t^{(1,5)}+\tilde f_2\\
f_3(x(t),y(t),z(t))&=(b-2c)t^{(1,6)}+\tilde f_3,
\end{align}
where $\tilde f_i$ stands for higher order terms with respect to the
$t$-adic valuation. Choose $b_\alpha, b_\beta, c_\alpha,$ $c_\beta$ so that
none of $f_1,f_2,f_3$ change sign between $\alpha$ and $\beta$. The
smallest $\nu_\alpha$-value of an element which changes sign between
$\alpha$ and $\beta$ is
\[
(1,4)+(0,4)=(1,5)+(0,3)=(1,8),
\]
so $\mu_\alpha=(1,8)$. 

Thus we have $f_i\notin<\alpha,\beta>$, but $f_i\in<\alpha,\beta>R[x,y,z]_{\mathfrak p}$.
\medskip

In this way we are naturally led to formulate a stronger version of CP, one 
which has exactly the same conclusion but with somewhat weakened hypotheses.

\begin{deft}\label{SCP}\textnormal{\textbf{(Strong Connectedness Property)}}
Let $\Sigma$ be a ring and
\[
\alpha,\beta \in \sper\ \Sigma
\]
two points having a common specialization $\xi$.
We say that $\Sigma$ has \textbf{the Strong Connectedness
Property at} $\alpha, \beta$ if given any $f_1,\ldots,f_s \in \Sigma \setminus
(\mathfrak{p}_\alpha \cup \mathfrak{p}_\beta)$ such
  that for all $i \in \{1,\ldots,s\}$,
\begin{equation} 
 \nu_\alpha(f_i) \leq
  \mu_\alpha,\ \ \nu_\beta(f_i) \leq
  \mu_\beta
\end{equation}
 and such that no $f_i$ changes sign between $\alpha$ and $\beta$, the points
$\alpha$ and $\beta$ belong to the same connected component of $\sper\
\Sigma \setminus \{f_1 \cdots f_s =0\}$.

We say that $\Sigma$ has the \textbf{Strong Connectedness Property} if it has the 
Strong Connectedness Property at $\alpha, \beta$ for all $\alpha, \beta \in 
\sper\ A$ having a common specialization. 

Let $n\in \N\setminus\{0\}$. We say that $\Sigma$ has \textbf{the Strong Connectedness Property up 
to height} $n$ if it has the Strong Connectedness property at $\alpha, \beta$ 
for all $\alpha, \beta \in \sper\ A$ having a common specialization such that 
$ht(<\alpha,\beta>) \le n$.
\end{deft}

\noi Let $\mathcal{X}$ be a class of rings. We use the following notation:

SCP($\mathcal{X}$) is the statement that every ring $\Sigma \in \mathcal{X}$ has the Strong  Connectedness property;

SCP$_n(\mathcal{X})$ is the statement that every ring in $\mathcal{X}$ of dimension at most $n$
has the Strong Connectedness property;

SCP$_{\le n}(\mathcal{X})$ is the statement that every ring in $\mathcal{X}$ has the Strong Connectedness property up to height $n$.

\begin{rek} \label{rek:local}
 One advantage of the Strong Connectedness Property is that its hypotheses behave well under localization at prime valuation ideals. Namely, let $\Sigma$, $\alpha$, $\beta$ and $f_1,\dots,f_s$ be as in Definition \ref{SCP}. Let $\xi$ be a common specialization of
$\alpha$ and $\beta$. Let $\alpha_0$ be the preimage of $\alpha$ under the natural inclusion $\sigma : \sper\
\Sigma_{{\mathfrak p}_\xi} \hookrightarrow \sper\ \Sigma$ and similarly for $\beta_0$. Then, for all $i \in \{1,\ldots,s\}$, 
\begin{align} 
 \nu_{\alpha_0}(f_i) &\leq
  \nu_{\alpha_0}(<\alpha,\beta>)= \nu_{\alpha_0}(<\alpha_0,\beta_0>),\\
\nu_{\beta_0}(f_i) &\leq \nu_{\beta_0}(<\alpha,\beta>)= \nu_{\beta_0}(<\alpha_0,\beta_0>)
\end{align}
and $f_i$ does not change sign between $\alpha_0$ and $\beta_0$. 
\end{rek}

\begin{conj} \textnormal{\textbf{(The Strong Connectedness Conjecture)}} \label{strong}

Let $\mathcal{X}$ be the class of all the regular rings. Then SCP($\mathcal{X})$ holds.
\end{conj}

Let $\Sigma$ be a ring and $\alpha, \beta \in \sper\ \Sigma$. Let the notation be as in Remark \ref{rek:local}. 
We have the following result (see \cite{LMSS3}):

\begin{thm} \label{thm:localsdcc}
If $\Sigma_{{\mathfrak{p}_\xi}}$ has the Strong Connectedness property at
$\alpha_0,\beta_0$, then $\Sigma$ has the Strong Connectedness property at
$\alpha, \beta$.
\end{thm}
\begin{cor} Let  $\mathcal{X}$ be a class of rings closed under localization.

\noi We have 
SCP$_n(\mathcal{X}) \implies$ SCP$_{\le n}(\mathcal{X})$ (the implication SCP$_{\le n}(\mathcal{X})\implies$
SCP$_n (\mathcal{X})$ is trivial and does not depend on $\mathcal{X}$ being closed under localization). 
\end{cor}

The Connectedness Property with $s=1$ will be referred to as \textbf{the Separation Property}. The \textbf{Separation Conjecture} asserts that any regular ring has the Separation Property; this is equivalent to the Connectedness Conjecture for $s=1$ and to the Pierce-Birkhoff Conjecture for $s=2$. Analogously to CP, we have the following stronger version of the Separation Property. 

\begin{deft} Let $\Sigma$ be a ring, $f$ a non-zero element of $\Sigma$ and $\alpha,\beta$ two points of $\sper\ \Sigma$ having a common specialization. Consider the following conditions:

(1) $\nu_\alpha(f) \leq \mu_\alpha$, $\nu_\beta(f) \leq \mu_\beta$ and $f$ does not change sign between $\alpha$ and $\beta$

(2) $\alpha$ and $\beta$ lie in the same connected component of $\{f\neq 0\}$.
\medskip

The \textbf{strong separation property} for the triple $(f,\alpha,\beta)$ is the implication (1)$\implies$(2). The ring
$\Sigma$ has the \textbf{strong separation property} if the strong separation property holds for any triple $(f,\alpha,\beta)$ as above.
\end{deft}

\begin{rek} \label{rek:centrecom}
	Assume that $\alpha$ and $\beta$ have a common specialization. Let $\xi$ be the most general such specialization. If $f(\xi) \neq 0$, then $\alpha$ and $\beta$ lie in the same connected component of $\{f\neq 0\}$, so $(f,\alpha,\beta)$ trivially satisfy the Strong Separation Property in this case. In the rest of the paper, we will tacitly assume that $f(\xi)=0$.  
\end{rek}

\noi Let $\mathcal{X}$ be a class of rings.  We use the following notation:

SP($\mathcal{X}$) is the statement that the Separation Property holds for all the rings in $\mathcal{X}$;
 
SSP($\mathcal{X}$) is the statement that the Strong Separation Property holds for all the rings in $\mathcal{X}$;

SP$_n(\mathcal{X})$ (resp. SSP$_n(\mathcal{X})$) is the statement that the Separation Property (resp. Strong Separation Property) holds for all the rings in $\mathcal{X}$ of dimension at most $n$.
 
SP$_{\le n}(\mathcal{X})$ (resp. SSP$_{\le n}(\mathcal{X})$) is the statement that the Separation Property (resp. the Strong Separation Property) holds for all rings $\Sigma \in \mathcal{X}$ and all $\alpha, \beta \in \sper\ \Sigma$ having a common specialization, such that $ht(\mathfrak{p}) \le n$.
\medskip

\begin{conj} \textnormal{\textbf{(Strong Separation Conjecture)}} Let $\mathcal{X}$ be the class of all regular rings. Then SSP($\mathcal{X}$) holds. 
\end{conj}
\medskip

As explained above, we have so far proved the following implications for any class of rings $\mathcal{X}$ and any natural number $n$:

\begin{equation}
\xymatrix{SCP_n(\mathcal{X}) \ar@{=>}[rr] \ar@{<=}[d]& & SSP_n (\mathcal{X}) \ar@{<=}[d] \\ SCP_{\le n}(\mathcal{X}) 
\ar@{=>}[rr]  \ar@{=>}[d] & &SSP_{\le n}(\mathcal{X}) \ar@{=>}[d] \\ CP_{\le n}(\mathcal{X}) \ar@{=>}[r] 
\ar@{=>}[d] & PB_{\le n}(\mathcal{X}) \ar@{=>}[r] \ar@{=>}[d] & SP_{\le n}(\mathcal{X}) \ar@{=>}[d] \\ 
CP_n(\mathcal{X}) \ar@{=>}[r] & PB_n(\mathcal{X}) \ar@{=>}[r] & SP_n(\mathcal{X})}
\end{equation}

\noi If, in addition, the class $\mathcal{X}$ is closed under localization at prime ideals, then the two upper vertical arrows of the above diagram are, in fact, equivalences. In other words, $SCP_n(\mathcal{X}) \iff SCP_{\le n}(\mathcal{X})$ and $SSC_n(\mathcal{X}) \iff SSC_{\le n}(\mathcal{X})$.
\medskip

\noi Consider the Euclidean space $R^{n+1}$ with coordinates $(x,z)$, where $x=(x_1,\ldots,x_n)$. Let $\pi:R^{n+1}\rightarrow R^n$ denote the natural projection onto the $x$-space. Let $D \subset R^n$ be a connected
semi-algebraic subset of $R^n$.  A \textbf{cylinder} in $R^{n+1}$ is a set of the form $C=\pi^{-1}(D)$ for some $D$ as above. Using the same notation as in \cite{BCR}, given a basic semi-algebraic set $F\subset R^{n+1}$, we denote by $\tilde{F}$ the subset of
$\sper\ B[z]$ defined by the same equations and inequalities as $F$ in $R^{n+1}$. More generally, if $F$ is a boolean combination of basic semi-algebraic subsets of $R^{n+1}$, $\tilde{F}$ is defined in the obvious way. Similarly, we denote by $\tilde{\pi}:\sper\ B[z]\rightarrow\sper\ B$ the natural projection corresponding to $\pi$. We will refer to $\tilde C$ as the \textbf{cylinder} lying over $\tilde
D$.

\begin{deft} \label{def:realbranch}
Let $g \in B[z]$. Let $\phi:D \to \sur{R}$ be a semi-algebraic continuous function and let $h=z-\phi : C \to \sur{R}$. 
If, for each $a \in D$, the element $\phi(a)$ is a simple root of the polynomial equation $g(a,z)=0$, we call $h$ a \textbf{branch} of $g$ over $D$. We call $h=z-\phi$ a \textbf{branch} over $D$ if $h$ is a branch of $g$ over $D$ for some $g$. It is called a \textbf{real branch} over $D$ if the image of $\phi$ is contained in $R$.
\end{deft}

Let $K$ be an algebraically closed field and $\nu$ a valuation on $K[z]$.
Let $g \in K[z]$ be a monic polynomial and let $g=\prod\limits_{i=1}^dg_i$ be the factorization of $g$ into linear factors.
\begin{deft} For $i\in\{1,\dots,d\}$, we say that $g_i$ is a $\nu$\textbf{-privileged factor} if $\nu(g_i) \geq \nu(g_j)$ for all
$j\in\{1,\ldots,d\}$. 
\end{deft}
\begin{deft} Assume given a real closed field $L\subset K$ such that $K=\sur{L}$. For $i\in\{1,\dots,d\}$, if $g_i \in L[z]$, we call $g_i$ \textbf{a real factor of} $g$. 
\end{deft}
Let $\gamma_0 \in \sper\ B$. Fix a real closure $\sur{B(\gamma_0)}_r$ of $B(\gamma_0)$ and an algebraic closure
$\sur{B(\gamma_0)}$ of $\sur{B(\gamma_0)}_r$, once and for all. Take $g = z^d +a_{d-1}z^{d-1} \cdots+ a_0 \in B[z]$. Let
$g=\prod\limits_{i=1}^d g_i$ be the factorization of $g$ into linear factors over $\sur{B(\gamma_0)}$. 

\begin{deft} We refer to the $g_i$ as $\gamma_0$-\textbf{branches}. If
	\[
	g_i \in \sur{B(\gamma_0)}_r[z],
	\]
	we say that $g_i$ is a \textbf{real} $\gamma_0$-\textbf{branch} of $g$.
\end{deft}

Take an element $\gamma \in \pi^{-1}(\gamma_0) \in \sper\ A$. A $\gamma$\textbf{-privileged} branch of $g$ is a
$\nu_\gamma$-privileged $\gamma_0$-branch of $g$. 
\medskip

In \S\ref{Value}, we will recall results from \cite{BCR} which canonically associate to each real branch $g_i$ of $g$ an element $g_i(\gamma) \in \sur{A(\gamma)}_r$.
\medskip

Consider a monic polynomial $f = z^d +a_{d-1}z^{d-1} \cdots+ a_0 \in A$. For a point $b\in D$, denote by $f(b)$ the polynomial
 $f(b) = z^d +a_{d-1}(b)z^{d-1} \cdots+ a_0(b) \in R[z]$.
\medskip

\noi\textbf{Notation.} For a polynomial $g\in A$ we will use the notation $g^{(k)}:=\frac1{k!}\frac{\partial g^k}{\partial z^k}$.

\begin{deft} \label{def:good} Assume that $\mathfrak{p}$ is a maximal ideal of $A$ so that $\mathfrak{p} = \text{supp}(\xi)$. Assume that there exist connected semi-algebraic sets $D \subset R^n$ and
$C=\pi^{-1}(D)\subset R^{n+1}$ as above, having the following properties:

(1) $\alpha, \beta \in \tilde{C}$

(2) For each $k\in\{0,\dots,d-1\}$, the number of real roots of $f^{(k)}(b)$, counted with or without multiplicity, is independent of the point $b\in D$.

In this situation we say that $\alpha,\beta$ \textbf{are in good position with respect to} $f,x,z$.
\end{deft}

\begin{rek}
Condition (2) of Definition \ref{def:good} is satisfied in the following situation. Assume that there are $d$ pointwise distinct continuous functions $\phi_j:D \to R$ such that $f = \prod\limits_{j=1}^d(z-\phi_j)$ in $C$ (in other words, $f$ has $d$ pointwise distinct real roots in $C$). Then, for all $i >0$, $\Delta\left(f^{(i)}\right) \neq 0$ and $f^{(i)}$ has $d-i$ pointwise distinct continuous real roots in $C$. In this case, any $\alpha, \beta \in C$ having a common specialization are in good position with respect to $f,x,z$.
\end{rek}

\begin{rek}
	The condition that $\mathfrak{p}$ is a maximal ideal of $A$ is equivalent to saying that $\alpha$ and $\beta$ have a unique common specialization whose support is a maximal ideal of $A$.
\end{rek}

We can now state the main theorem of this paper.
 
\begin{thm} \label{thm:main}
 Assume that $\alpha,\beta$ are in good position with respect to $f,x,z$. Then 
the Strong Separation Property holds for $f$, $\alpha$ and $\beta$. 
\end{thm}

At the end of the paper, we will explain that the hypothesis of good position in the theorem can be relaxed somewhat (see Remark \ref{rek:extension}).
\medskip
 
\noi \textbf{Open question.} We do not know whether, given $f \in A$ and $\alpha, \beta \in \sper\ A$, there exists an automorphism $\sigma : A \to A$ such that, letting $\tilde{x}_j =\sigma(x_j)$ and $\tilde{z}=\sigma(z)$, the points $\alpha$ and $\beta$ are in good position with respect to $f, \tilde{x}$ and $\tilde{z}$. 
\medskip

\begin{rek} We note that the following slightly more general result holds.

Let $B_D \subset R(x)$ denote the ring of rational functions having no poles 
in $\sur{D}$.

To each polynomial $g$ in $A$ whose first coefficient has no zeroes in $\sur{D}$ we can naturally associate a monic polynomials in
$B_D[z]$, namely, $g$ divided by its leading coefficient. 

The above definitions, in particular, Definition \ref{def:good}, extend in an obvious way to polynomials in $B_D[z]$. Theorem \ref{thm:main} holds for $f \in B_D[z]$, alternatively, for non-monic $f \in A$ whose first coefficient has no zeroes in $\sur{D}$.  The proof is exactly the same as the proof of Theorem \ref{thm:main} given in the present paper. We chose to work in the more restrictive setup of monic polynomials in $A$ in order not to overburden the notation.
\end{rek}

\begin{rek}
	As we explained earlier the reason for introducing strong versions of all the conjectures is the fact that the strong versions are stable under the localization at $\mathfrak{p}$. This would allow us to proceed by induction on the dimension of $A$ and reduce the case when $\mathfrak{p}$ is not maximal to the case when $\mathfrak{p}$ is maximal by localization. The difficulty is that in the preset paper we use in an essential way the hypothesis that the ground field is real closed. Localization at $\mathfrak{p}$ destroys this hypothesis. If the results of the present paper could be generalized to a non real closed ground field $R$, the hypothesis on the maximality of $\mathfrak{p}$ would become unnecessary.   
\end{rek}

\begin{rek} In the proof of Theorem \ref{thm:main}, we may assume that the polynomial $f$ is reduced. Indeed, suppose the 
Theorem is true when $f$ is reduced. Assume that $\tilde{f} \in A$ 
satisfies $\nu_\alpha(\tilde{f}) \leq \mu_\alpha$, $\nu_\beta(\tilde{f}) 
\leq\mu_\beta$ and that $\tilde{f}$ does not change sign between $\alpha$ and $\beta$.
Let $f = \tilde{f}_{red}$. 
We have $\nu_\alpha(f)  \le \nu_\alpha(\tilde{f}) \le \mu_\alpha$ 
and the same for $\beta$.
If
\begin{equation}
\nu_\alpha(f) = \nu_\alpha(\tilde{f}),\label{eq:vlaueequality}
\end{equation}
then $\nu_\alpha(\tilde{f}/f) = 0$, so $\tilde{f}/f$ is a unit in $R_\alpha$ and $R_\beta$. Therefore $f$ does not change sign between $\alpha$ and $\beta$. Thus, in all cases (that is, regardless of whether the equality (\ref{eq:vlaueequality}) holds) $f$ satisfies the hypotheses of the Strong Separation Conjecture. Then $\alpha$ and $\beta$ belong to the same connected component of $\{f\ne0\}$ and hence also to the same connected component of $\{\tilde f\ne0\}$, as desired.
\end{rek}

This paper is organized as follows.

In \S 2 we let $K$ be an algebraically closed field. Let $g = \sum\limits_{i=0}^d\ a_iz^i$, where $a_i \in K$ and $a_d=1$, be a monic polynomial. Let $\nu$ be a valuation of $K[z]$.  
\medskip

We define the notion of Newton polygon of $g$. The main result of \S 2, Proposition 
\ref{prop:privbranch}, says that if $k\in\{1,\dots,d-1\}$ and $\nu_\gamma\left(g^{(i)}\right) \ge \nu_\gamma(g)$ for all $i\le k$ 
then, for each $\gamma_0$-branch $g_j$ of $g$ and each $\gamma$-privileged branch $h$ of 
$g^{(k)}$, we have $\nu_\gamma(h) > \nu_\gamma(g_j)$.  
\medskip

In \S\ref{semialgebraicity} we start with a monic polynomial $g$ in the variable $z$ whose coefficients are continuous functions over a semi-algebraic set $U\subset R^n$. We assume that for $a\in U$ the number $s$ of real roots of $g(a)$, counted with multiplicities, is independent of $a$. We define continuous functions $\phi_i:U\rightarrow R$, $1\le i\le s$, such that for all $a\in U$ we have
$g(a)=\prod\limits_{i=1}^s(z-\phi_i(a))$. This result is well known in the case when $R=\mathbb R$, but we did not find it in the literature in the case of an arbitrary real closed field $R$. If, in addition, the coefficients of $g$ are semi-algebraic functions on $U$ then the $\phi_j$ can also be chosen to be semi-algebraic : \cite{DK} Lemma 1.1. 

In \S\ref{Value} we recall results from \cite{BCR} which canonically associate to each real branch $g_i$ of $g$ an element $g_i(\gamma) \in \sur{A(\gamma)}_r$.
\medskip

In \S\ref{Sec:reel}  we study, for a point $\gamma\in\sper\ A$, the extension $\nu_{\bar\gamma}$ of $\nu_\gamma$ to
$\sur{A(\gamma)}$. The main result of \S\ref{Sec:reel} is Proposition \ref{lem:reim}; it says that this extension is unique and that for every branch $h$ over $D$ we have
\begin{equation} 
 \nu_{\sur{\gamma}}(h) = \min\{\nu_{\sur{\gamma}}(Re\, h), 
\nu_{\sur{\gamma}}(Im\, h) \}.
\end{equation}
Most of \S\ref{Sec:Rolle} is devoted to using the Newton polygon to prove comparison 
results between quantities of the form $\nu_\gamma(g_j)$ and $\nu_\gamma(h)$ 
where $g_j$ is a $\gamma_0$-branch of $g$ and $h$ is an $\gamma_0$-branch of $g'$, as well as relating 
inequalities of size to inequalities of values. Lemma \ref{lem:conv2} says that 
if $h_1,h_2$ are two real branches such that $0< h_1(\gamma) < h_2(\gamma)$ then 
$\nu_\gamma(h_1(\gamma)) \ge \nu_\gamma(h_2(\gamma))$. 

If $g_1,g_2$ are two real branches of $g$ and $h$  a real branch of $h'$ between 
$g_1$ and $g_2$, then $\gamma$ cannot have strictly higher contact with both 
$g_1$ and $g_2$ than it does with $h$. The equidistance Lemma (Lemma 
\ref{lem:equid3}) is a valuation-theoretic generalization of this fact to the 
case when the branches are not necessarily real. 

\S\ref{Sec:Reduc} is devoted to  reducing the problem to the case when $\Gamma_\alpha \cong \Z$ and
$\Gamma_\beta \cong \Z$ and $k_\alpha = k_\beta=R$.

Finally, in \S\ref{Sec:Proof} we use the results of the preceding sections to complete the proof of 
Theorem \ref{thm:main} by induction on $deg(f)$.
\medskip

\noi\textbf{Acknowledgement.} We thank Michel Coste for sharing his insights on semi-algebraicity of roots of polynomials over real closed fields.

We would also like to thank the Institute of Mathematics of UNAM in Cuernavaca (Mexico) for inviting D. Schaub in February 2018; the work on this paper was completed during our stay there.  

\section{Graded algebra and Newton polygon}
 
 Let $K$ be an algebraically closed field and $z$ an independent variable. Fix a valuation $\nu$ of $K(z)$ with value group $\Gamma$.
 


\noi Let $\nu_z$ denote the valuation of $K(z)$ defined by 
\[\nu_z\left(\sum\limits_i 
b_iz^i\right) = \min\left\{\left.\nu\left(b_iz^i\right)\ \right|\ b_i\neq 0\right\}.\]
For each $g \in K[z]$, we have $\nu_z(g) \leq \nu(g)$. We will write 
$\init_z$ instead of $\init_{\nu_z}$.

\begin{rek} \label{rek:isom}
	Let $X$ be an independent variable. Since $\left.\nu\right|_K = \left.\nu_z\right|_K$ and $\init_zz$ is transcendental over
$\gr_zK$ in $\gr_z(K[z])$, we have a natural isomorphism 
\[
(\gr_\nu K)[X] \cong \gr_z(K[z])
\]
of graded algebras.
\end{rek}
	
For $g=\sum\limits_{j=0}^da_jz^j$, let $S_z(g) = \left\{i\in\{0,\dots,d\}\ \left|\ 
\nu\left(a_iz^i\right) = \nu_z(g)\right.\right\}$. We have
\[
\init_zg =\sum_{i \in S_z(g)} \init_z(a_i) \init_z(z)^i.
\]
Assume that there exists $g \in K[z]$ such that 
\begin{equation} \label{eq:inegal}
 \nu_z(g) < \nu(g).
\end{equation}
 For such a $g$, the polynomial $\init_zg$ is not a monomial.
\smallskip

Let $\delta(g,z)=\deg(\init_zg)$ and $\delta = \delta(g,z)$.

\noi The inequality (\ref{eq:inegal}) is equivalent to saying that
\[
\sum_{i \in S_z(g)} \init_\nu(a_i)\init_\nu(z)^i=0.
\]
\noi Since $K$ is algebraically closed, the polynomial $\bar{g}(X) := \sum\limits_{i \in S_z(g)} \init_\nu a_i X^i$ factors into linear factors in $\gr_\nu K[X]$ :
\[
\bar{g}(X) = \init_\nu a_\delta \prod_{j=1}^\delta (X-\psi_j).
\]
Thus there exists a unique $\bar{\phi} \in \gr_\nu K$ such that $\init_\nu z = \bar{\phi}$. Take a representative $\phi \in K$ such that
$\bar{\phi} = \init_\nu \phi$. We have $\nu(z-\phi) > \nu(z)$. 

We summarize the above considerations in the following Remark:
\begin{rek} \label{rek:algebric}
 Take $g \in K[z]$. We have the following 
implications: 
\begin{align*} \nu_z(g) < \nu(g) &\implies \init_\nu(z) \text{ is algebraic 
over } \gr_\nu(K) \iff \init_\nu(z)  \in \gr_\nu(K) \\ &\iff \text{there exists }\phi \in 
K \text{ such that } \nu(\phi) =\nu(z) < \nu(z-\phi).
\end{align*}
The second implication uses the fact that $K$ is algebraically closed; the first and the third one are valid without any hypotheses on $K$. 
\end{rek}

Assume that the strict inequality (\ref{eq:inegal}) holds. Consider the change of 
variables $\tilde{z}=z-\phi$ as above.  Write $g= \sum\limits_i 
\tilde{a}_i\tilde{z}^i$. Let $S_{z,\tilde{z}}(g) = \left\{i\ 
\left|\ \nu_z\left(\tilde{a}_i\tilde{z}^i\right)= 
\nu_z(g)\right. \right\}$. Let
\[
\kappa = \max S_{z,\tilde{z}}(g).
\]
We have $\nu_z\left(\tilde{a}_i\tilde{z}^i\right) = \nu\left(\tilde{a}_i\right)+ i \nu(z)$ 
because $ \nu_z\left(\tilde{a}_i\right) = \nu\left(\tilde{a}_i\right)$ and
\[
\nu_z\left(\tilde{z}\right) = \nu_z(z-\phi) = \min\{\nu(z),\nu(\phi)\} = \nu(z) = \nu(\phi).
\]
\begin{lem} \label{lem:ztildez} We have:

 (1) $\nu_z(g) = \min \left\{\left.\nu_z\left(\tilde{a}_i\tilde{z}^i\right)\ \right|\ \tilde{a}_i 
\neq 0\right\}$.

(2) $\init_zg = \sum\limits_{i \in S_{z,\tilde{z}}(g)} \init_z\left(\tilde{a}_i\right)\init_z(z-\phi)^i$. 
\end{lem}

\noi\textit{Proof:}  (1) The fact that $\nu_z(g) \geq \min\left\{\left.\nu_z\left(\tilde{a}_i\tilde{z}^i\right)\ \right|\
\tilde{a}_i\neq0\right\}$ follows from the definition of valuation. Now, $\nu_z(g) =\min\left\{\left.\nu\left(a_iz^i\right)\ \right|\ 
a_i\neq0\right\}$. Replacing $z$ by $\tilde{z}+ \phi$ in $g$ and expanding in $z$, we see that $\tilde{a}_i$ is a sum of terms of the form 
$ca_{i+k}\phi^k$ where $c \in \N$. Hence
$\nu\left(\tilde{a}_i\tilde{z}^i\right)\geq\min\limits_{k\in\N}\{\nu(a_{i+k}\phi^k z^i)\}\geq
\min\limits_{0\le j\le d}\left\{\nu\left(a_jz^j\right)\right\} = \nu_z(g)$. This proves (1). 

\noi From (1), we deduce that $S_{z,\tilde{z}}(g) = \left\{i\ \left|\ \nu_z\left(\tilde{a}_i \tilde{z}^i\right)= \nu_z(g)\right.\right\}$. 
We have
\[
g = \sum_{i \in  S_{z,\tilde{z}}(g)}\tilde{a}_i \tilde{z}^i + h
\]
with $\nu_z(h) > \nu_z(g)$. This proves (2).  $\Box$

\begin{rek}
 With the above notation, we have $\kappa=\delta$. 
\end{rek}

\begin{lem} \label{lem:tilde}
 Consider two integers $i,j \in \{0,\ldots,d\}$. Assume that 
\begin{equation} \label{eq:inegal2} \nu\left(\tilde{a}_i\tilde{z}^i\right) \leq \nu\left(\tilde{a}_j \tilde{z}^j\right) \text{ and } \nu_z\left(\tilde{a}_i\tilde{z}^i\right) \geq \nu_z\left(\tilde{a}_j \tilde{z}^j\right).  \end{equation}
We have:

(1) $i \leq j$;

(2) If at least one of the inequalities (\ref{eq:inegal2}) is strict, then $i<j$. 
\end{lem}

\noi Proof : We have $\nu\left(\tilde{a}_i)+ i\nu(\tilde{z}\right) =\nu\left(\tilde{a}_i\tilde{z}^i \right)  \leq
\nu\left(\tilde{a}_j \tilde{z}^j\right) =\nu\left(\tilde{a}_j\right) + j\nu\left(\tilde{z}\right)$, so
\[
\nu\left(\tilde{a}_i\right) - \nu\left(\tilde{a}_j\right) \leq (j-i) \nu\left(\tilde{z}\right).
\]
On the other hand, 
\[
\nu\left(\tilde{a}_i\right) + i \nu(z)  = \nu\left(\tilde{a}_i\right)+ i\nu_z\left(\tilde{z}\right) =\nu_z\left(\tilde{a}_i \tilde{z}^i\right)\geq 
\nu_z\left(\tilde{a}_j \tilde{z}^j\right) = \nu\left(\tilde{a}_j\right) + j \nu(z),
\]
so $(j-i) \nu(z) \leq \nu\left(\tilde{a}_i\right)-\nu(\tilde{a}_j)$. We obtain
\[
(j-i) \nu(z) \leq \nu\left(\tilde{a}_i\right)-\nu\left(\tilde{a}_j\right) \leq  (j-i)\nu\left(\tilde{z}\right).
\]
From this we deduce both (1) and (2).  $\Box$
\medskip

Let $\tilde{\delta} = \delta(g,\tilde{z}) = \deg(\init_{\tilde{z}}(g))$. 

\begin{lem}\label{lem:delta} We have:

(1) $\tilde\delta \leq  \delta$;

(2) If equality holds in (1), then $\init_zg= \init_z(\tilde{a}_\delta) 
\init_z(z-\phi)^\delta $.
\end{lem}

\noi Proof : We have $\nu\left(\tilde{a}_{\tilde{\delta}}\tilde{z}^{\tilde{\delta}}\right) \leq
\nu\left(\tilde{a}_\delta \tilde{z}^\delta\right)$, because $\tilde{\delta} \in S_{\tilde{z}}(g)$, so
$\nu\left(\tilde{a}_{\tilde{\delta}}\tilde{z}^{\tilde{\delta}}\right)$  is minimal among all the $\nu\left(\tilde{a}_i\tilde{z}^i\right)$. 

On the other hand, we have $\nu\left(\tilde a_\delta\right)=\nu\left(a_\delta\right)$. By Lemma \ref{lem:ztildez},
\[
\nu_z\left(\tilde{a}_{\tilde{\delta}}\tilde{z}^{\tilde{\delta}}\right) \geq\nu_z(g)=\nu_z\left(a_\delta z^\delta\right)=
\nu_z\left(\tilde{a}_\delta \tilde{z}^\delta\right).
\]
Applying (1) of Lemma \ref{lem:tilde} with $\tilde{\delta}=i$ and $\delta=j$, we deduce (1). 

Now, if equality holds in (1), $S_{z,\tilde{z}}(g) = \{\delta\}$. Indeed, let $i \in 
S_{z,\tilde{z}}(g)$; then $i \le \delta$. Now, by definition of 
$\tilde{\delta}$, we have 
$\nu\left(\tilde{a}_i\tilde{z}^i\right) \ge \nu\left(\tilde{a}_{\tilde{\delta}} 
\tilde{z}^{\tilde{\delta}}\right) = \nu\left(\tilde{a}_\delta 
\tilde{z}^\delta\right)$ (as $\delta=\tilde{\delta}$). On the other hand, 
because $i,\delta \in S_{z,\tilde{z}}(g)$, we have 
$\nu_z\left(\tilde{a}_i\tilde{z}^i\right) = 
\nu_z\left(\tilde{a}_\delta \tilde{z}^\delta\right)$. Applying (1) of Lemma 
\ref{lem:tilde} to the pair $(i,\delta)$, we deduce that $\delta \le i$. Thus $\delta = i$ and $S_{z,\tilde{z}}(g) = \{\delta\}$. Now the conclusion follows from Lemma \ref{lem:ztildez} (2).
$\Box$ \medskip

\begin{rek} \label{rek:delta0} Let $g(z) = a_0+ a_1z+ \cdots+a_dz^d \in K[z]$.

(1) We have the following implications: $\delta(g,z) = 0 \iff \nu(a_0) < \nu(a_jz^j) $ for all $j>0$ 
$\implies \nu_z(g) = \nu(g)$.  

(2) Assume that $\init_\nu(z) \in 
\gr_\nu(K)$. Let $\phi \in K$ be as in the Remark \ref{rek:algebric} and put 
$z_1=z-\phi$. We have $\nu_z(g) = \nu(g) \implies \delta(g,z_1)=0$.
\end{rek}

\begin{lem}
 \label{lem:stationne}
 Let $\Sigma$ be a noetherian domain and $\mu$ a rank 1 valuation of the field of fractions of $\Sigma$, non-negative on $\Sigma$. Then every bounded subset of the semi-group $\mu(\Sigma\setminus\{0\})$ is finite.
\end{lem}

\noi\textit{Proof:} Localizing $\Sigma$ at the center of $\mu$ does not change the problem. Assume that $(\Sigma,\mathfrak{m})$ is local and $\mu$ is centered at $\mathfrak{m}$. Take a subset $T \subset \mu(\Sigma\setminus\{0\})$ such that $T \le \beta$ for some
$\beta \in \mu(\Sigma\setminus\{0\})$. Let $I_\beta = \{y \in\Sigma \ |\ \mu(y) \ge \beta\}$. Since rk$(\mu) =1$, we have
$\mathfrak{m}^n \subset I_\beta$ for some $n \in \N$. A chain of elements $\beta_1 < \beta_2 < \ldots \le \beta$ of $T$ induces a chain of submodules
\[\frac{I_{\beta_\ell}}{\mathfrak{m}^n}\subset\frac{I_{\beta_{\ell-1}}}{\mathfrak{m}^n}\subset\ldots\subset
\frac{I_{\beta_1}}{\mathfrak{m}^n} \subset \frac{\Sigma}{\mathfrak{m}^n}.\] Hence
\[\#T\le \text{length}\frac\Sigma{\mathfrak{m}^n}<\infty.\] $\Box$
\medskip

The set of isolated subgroups of $\Gamma$ is naturally ordered by inclusion. In the applications the rank of $\Gamma$ will be finite by Abhyankar's inequality. In particular, the set of isolated subgroups of $\Gamma$ will be well ordered.
\medskip

From now till the end of this section, assume that $char\ K=0$.

\begin{lem} \label{lem:gene} 
Assume that the set of isolated subgroups of $\Gamma$ is well-ordered. Fix a polynomial
\begin{equation}
g(z)=\sum\limits_{i=0}^d a_i z^i \in K[z].\label{eq:polynomialg}
\end{equation} 

(1) There exists $\phi \in K$ such that, letting $\tilde{z} = z-\phi$, we 
have $\nu_{\tilde{z}}(g) = \nu(g)$. 

(2) Let $\tilde z$ be as in (1). Assume that $\init_\nu(\tilde{z}) \in \gr_\nu(K)$. Then there exists $\phi^* \in K$ such that, letting
$z^* = \tilde{z}-\phi^*$, we have $\delta(g,z^*)=0$. 
\end{lem}

\noi\textit{Proof:} By Remark \ref{rek:delta0}, (1) implies (2). Let us prove (1). Let $\delta=\delta(g,z)$.

Let $\Lambda_0$ be the smallest isolated subgroup of $\Gamma$ such that $\nu(g) - \nu_z(g) \in \Lambda_0$. If $\Lambda_0 = (0)$, there is nothing to prove. Assume that $(0)\subsetneqq\Lambda_0$. Let $\Lambda_-$ be the union of all the proper isolated subgroups of $\Lambda_0$. We have $\Lambda_-\subsetneqq \Lambda_0$ since
\[
\nu(g) - \nu_z(g) \in \Lambda_0 \setminus \Lambda_-.
\]
Thus $\Lambda_-$ is the greatest proper isolated subgroup of $\Lambda_0$.   

It is sufficient to show that there exists $\phi \in K$ such that, letting 
$\tilde{z}=z-\phi$, we have
\begin{equation}
\nu(g) - \nu_{\tilde{z}}(g) \in \Lambda_-.         
\end{equation}
The proof will then be finished by transfinite induction on $\Lambda_0$. 

Let $R_\nu \subset K(z)$ denote the valuation ring of $\nu$. Let
\begin{align}
\mathbf{P}_{\Lambda_-} &= \{y \in R_\nu \ |\ \nu(y) \notin \Lambda_-\} \\
\mathbf{P}_{\Lambda_0} &= \{y \in R_\nu \ |\ \nu(y) \notin \Lambda_0\}.
\end{align}
Write $\nu =\nu_0\circ\nu_- \circ \mu$ where $\nu_0$ is the valuation of $K(z)$ with 
valuation ring $(R_\nu)_{\mathbf{P}_{\Lambda_0}}$, $\nu_-$ is the rank one valuation of the residue field $\kappa(\mathbf{P}_{\Lambda_0})$ of $\mathbf{P}_{\Lambda_0}$ with 
valuation ring \[ \frac{(R_\nu)_{\mathbf{P}_{\Lambda_-}}}{\mathbf{P}_{\Lambda_0}(R_\nu)_{\mathbf{P}_{\Lambda_-}}}\] and $\mu$ is the valuation of the residue field  $\kappa(\mathbf{P}_{\Lambda_-})$ of $\mathbf{P}_{\Lambda_-}$  with valuation ring \[\frac{R_\nu}{\mathbf{P}_{\Lambda_-}}.\]
Replacing $\nu$ by $\nu_0\circ\nu_-$ does not change the problem. In this way, we may 
assume that $\Lambda_-=(0)$ and $rk(\Lambda_0)=1$. 

In what follows, we will consider changes of variables of the form 
\begin{equation} 
\label{eq:transfo} \tilde{z} = z-\phi,\end{equation} such that 
$\nu_{\tilde{z}}(g) \ge \nu_z(g)$. 

We proceed by induction on $\delta(g,z)$. The case $\delta=0$ is given by Remark \ref{rek:delta0}. Assume that $\delta >0$.  If
\begin{equation} 
\nu_z(g) = \nu(g), 
\end{equation}
that is, the conclusion of Lemma \ref{lem:gene} (1) holds with $\phi=0$ and 
$\tilde{z}=z$, there is nothing to prove. Assume 
 \begin{equation}\label{ineg:nu}
\nu_z(g) < \nu(g).
\end{equation}
\noi By the algebraic closedness of $K$, $\init_z (g)$ decomposes into 
linear factors in
\[
\init_z K[\init_z z] = \gr_z K[z].
\]
Hence, by (\ref{ineg:nu}) and Remark \ref{rek:algebric}, there exists $b\in K$ such 
that \begin{equation} 
\nu(z)=\nu(b) < \nu(z-b).
\end{equation}
\smallskip

\noi\textbf{Claim 1.} It is sufficient to prove (1) in the case when
\begin{equation}
\nu(z)=0\label{eqnuz=0}
\end{equation}
and
\begin{equation}
\nu(a_i)\ge\nu(a_\delta)\quad\text{ for all }i\in\{0,\dots,d\}.\label{eq:nuaige0}
\end{equation}
\noi\textit{Proof of Claim 1.} Assume that (1) of the Lemma is known in the case when (\ref{eqnuz=0}) and (\ref{eq:nuaige0}) hold. Put $z_b:=\frac zb$. We have
\begin{equation}
\nu(z_b)=0.\label{eq:nuzb>0}
\end{equation}
For $i\in\{0,\dots,d\}$, put $a_{ib}:=b^ia_i$. Using the definition of $\delta$, for each $i\in\{0,\dots,d\}$ we obtain
\begin{equation}
\nu\left(a_{ib}z_b^i\right)=\nu\left(a_iz^i\right)\ge\nu\left(a_\delta z^\delta\right)=\nu\left(a_{\delta b}z_b^{\delta}\right).\label{eq:ineqmonom}
\end{equation}
Together with (\ref{eq:nuzb>0}), this implies that $\nu\left(a_{ib}\right)\ge\nu\left(a_{\delta b}\right)$ for all $i\in\{0,\dots,d\}$.

Let $g_b:=\sum\limits_{i=0}^da_{ib}z_b^i$.  By assumption and in view of (\ref{eq:nuzb>0}), there exists $\phi_b\in K$ such that, letting $\tilde z_b=z_b-\phi_b$, we have
\begin{equation}
\nu_{\tilde z_b}(g_b)=\nu(g_b).\label{eq:nuzb=nugb}
\end{equation}
Put $\phi:=b\phi_b$. We have $\init_{\tilde z_b}g_b=(\init_{\tilde z}g)_b$. In particular, $S_{\tilde z}(g)=S_{\tilde z_b}(g_b)$. Now, (\ref{eq:nuzb=nugb}) is equivalent to saying that $\sum\limits_{i\in S_{\tilde z_b}(g_b)}\init_\nu\left(\tilde a_{ib}\tilde
z_b^i\right)\ne0$. Then $\sum\limits_{i\in S_{\tilde z}(g)}\init_\nu\left(\tilde a_i\tilde z^i\right)\ne0$, so that
\begin{equation}
\nu_{\tilde z}(g)=\nu(g).\label{eq:nuz=nug}
\end{equation}
This completes the proof of Claim 1.
\medskip

From now till the end of the proof of Lemma \ref{lem:gene}, assume that (\ref{eqnuz=0})--(\ref{eq:nuaige0}) hold. 
 
Replacing $g(z)$ by $\frac{g(z)}{a_\delta}$ does not change the problem. In this way, we may assume that
\begin{equation}
a_\delta=1\label{eq:adelta=1}
\end{equation}
and
\begin{equation}
\nu(a_i)\ge0\quad\text{ for all }i \in \{0,\ldots,\delta-1\}.\label{eq:nuaige0vraiment}
\end{equation}
\smallskip

Let \begin{align} 
P_{\Lambda_0} &:= \{w \in \Q[a_0,\ldots,a_d]\ |\ \nu(w) \notin \Lambda_0 \} = \mathbf{P}_{\Lambda_0} \cap \Q[a_0,\ldots,a_d] \\ 
\mathfrak{m}_0 &:= \{w \in \Q[a_0,\ldots,a_d]\ |\ \nu(w) >0 \} = \mathfrak{m}_\nu \cap \Q[a_0,\ldots,a_d].
 \end{align}
Assumption (\ref{eq:nuaige0vraiment}) imply that $P_{\Lambda_0}$ and $\mathfrak{m}_0$ are prime ideals of $\Q[a_0,\ldots,a_d]$. The valuation $\nu$ induces a valuation of $\Q(a_0,\ldots,a_d)$, centered at $\mathfrak{m}_0$.
\medskip

 Let $z_0:=z$. 
\smallskip

\noi Assume that, for a certain integer $s \ge 0$, we have constructed elements
\[
\phi_1,\phi_2,\ldots,\phi_{s-1}\in K
\]
and monic linear polynomials $z_0,\dots,z_s\in K[z]$, having the following properties (if $s=0$ we adopt the convention that both sets
$\{\phi_0,\dots,\phi_{s-1}\}$ and $\{z_1,\dots,z_{s-1}\}$ are empty, only $z_0$ is defined) :
 
 1) $z_{i+1} = z_i - \phi_i$, $i \le s-1$; we have
\begin{equation}\label{eq:comp1}
\nu(z_i) = \nu(\phi_i) < \nu(z_{i+1});
\end{equation}

 2) $g\in \Q[a_0,\ldots,a_d]_{\mathfrak{m}_0}[z_i]$ and $\phi_0,\dots,\phi_{s-1}\in \Q[a_0,\ldots,a_d]_{\mathfrak{m}_0}$;
 
 3) write $g= \sum\limits_{j=0}^da_{j,s} z_{s}^j$; we have
$a_{\delta,s}\in1+\mathfrak{m}_0\Q[a_0,\ldots,a_d]_{\mathfrak{m}_0}$.
\smallskip

In the case $s=0$ condition 3) clearly holds and 1) and 2) are vacuously true, since only $z_0$ is defined and the set $\{\phi_0,\dots,\phi_{s-1}\}$ is empty.
 
 \begin{rek} \label{rek:stats} 1. By (\ref{eq:comp1}) we have $\nu(\phi_0) < \nu(\phi_1) < \ldots < \nu(\phi_{s-1})$.
  
  2. By Lemma \ref{lem:delta}, $\delta(g,z) \ge \delta(g,z_1) \ge \ldots \ge \delta(g,z_s)$.
 \end{rek}

 If $s\ge1$ and $\delta(g,z_s) < \delta(g,z)$, the proof is finished by induction on $\delta(g,z)$. Assume that
\begin{equation}\label{eq:delta}
\delta(g,z_s) = \delta(g,z)= \delta;
\end{equation}
if $s>0$ this implies that $\delta(g,z_{s-1}) =\delta(g,z_s)$. 
 \smallskip
  
  If 
  \begin{equation} 
  \nu_{z_s}(g) = \nu(g), 
  \end{equation}
 that is, the conclusion of Lemma \ref{lem:gene} (1) holds with $\phi=\sum\limits_{j=1}^{s-1}\phi_j$ and 
$\tilde{z}=z_s$, there is nothing more to prove. Assume that
  \begin{equation}
\nu_{z_s}(g) < \nu(g).
\end{equation}
    
Let $X$ be an independent variable and consider the polynomial \begin{equation}\sur{g}(X) := \sum_{i \in S_{z_s}(g)} \init_\nu a_{i s}X^i. \end{equation} 

\noi By (\ref{eq:delta}) we have $\deg_X(\sur{g}) = \delta$. By Remark \ref{rek:algebric} and considerations which precede it, we have $\sur{g}(\init_\nu z_s) =0$. 

Since $K$ is algebraically closed , we can factor $\sur{g}$ into linear factors over $\gr_\nu K$. If $\sur{g}$ is not of the form \begin{equation}\label{eq:gbar}
\sur{g}= \init_\nu a_{\delta,s} (X- \psi)^\delta,
\end{equation}
take an element $\phi_s \in K$ such that $\init_\nu z = \init_\nu \phi_s$. Put $z_{s+1}:= z_s-\phi_s$. 

By Lemma \ref{lem:delta} (2) and in  view of Remark \ref{rek:isom}, we have
$\delta(g,z_{s+1}) < \delta$ and the proof is finished by induction on $\delta$. 

Assume that $\sur{g}$ is of the form (\ref{eq:gbar}). By Newton binomial theorem and in view of (\ref{eq:adelta=1}), equating the coefficients of $\init_{z_s} (z_s)^{\delta-1}$ on the right and left hand sides of (\ref{eq:gbar}), we see that 
\[
\psi = -\frac{\init_{z_s}(a_{\delta-1,s})}{\delta}.
\]
Define $\phi_s$ to be $-\frac{a_{\delta-1,s}}{\delta}$. By definitions, (\ref{eq:comp1}) holds for $i=s$. 

Repeat the procedure to construct a sequence $\phi_0, \phi_1,\phi_2, \ldots \in K$ such that
\begin{equation} \label{eq:suite}
      \nu(\phi_0) < \nu(\phi_1) < \nu(\phi_2) < \dots
     \end{equation} 
having the properties 1), 2), 3) preceding Remark \ref{rek:stats}.
\medskip

It remains to show that at some point of this construction we obtain
\[
\nu_{z_i}(g) = \nu(g)
\]
or
\[
\delta(g,z_i) < \delta(g,z).
\]
In both cases, the Lemma will be proved. Therefore, the proof of the Lemma is reduced to the following Claim:
\smallskip

\noi \textbf{Claim 2}: The above procedure cannot continue indefinitely. 
\smallskip

\noi Proof of the Claim 2: We give a proof by contradiction. Assume that the procedure does not stop, in other words, the sequence
$\{\phi_i\}$ is infinite. By construction,
\[
\phi_i \in \Q[a_0, \ldots, a_d]_{\mathfrak m_0}\quad\text{ for all }i \in \N.
\]
\noi Write $\nu_0=\theta\circ\epsilon$, where $\theta$ is a valuation of $\Q[a_0, \ldots, a_d]$, centered at $P_{\Lambda_0}$. We have $rk(\epsilon) = rk(\Lambda_0)=1$. Since 
 \[
\frac{\Q[a_0,\ldots,a_d]_{\mathfrak{m}_0}}{P_{\Lambda_0}\Q[a_0,\ldots,a_d]_{\mathfrak{m}_0}}
\]
is noetherian, (\ref{eq:suite}) and Lemma \ref{lem:stationne} imply that the sequence 
$\{\epsilon(\frac{\phi_i}{\phi_0})\}_{i\in \N}$ is unbounded in $\Lambda_0$. Hence so is the sequence
$\{\nu(\frac{\phi_i}{\phi_0})\}_{i\in \N}$.
\smallskip

Take $i$ sufficiently large so that 
\begin{equation} \label{eq:comp2}
 \delta \nu\left(\frac{\phi_i}{\phi_0}\right) > \nu(g) - \nu_z(g) .
 \end{equation} 

We have 
\begin{align}
 \nu_z(g) &= \delta \nu(\phi_0) \\ \nu_{z_i}(g) &= \delta \nu(\phi_i).
\end{align}

Hence \begin{equation}
       \nu(g) - \nu_z(g) < \delta(\nu(\phi_i)-\nu(\phi_0)) = \nu_{z_i}(g) - \nu_z(g) < \nu(g) - \nu_z(g),
      \end{equation}
which is a contradiction. This completes the proof of the Claim and with it the Lemma. $\Box$
\medskip

\begin{cor} \label{cor:gene2}
 Given a finite family of polynomials $f_1,\ldots, f_s  \in K[z]$, the following hold.

(1) There exists $\phi \in K$ such that letting $\tilde{z} = z - \phi$, we have
\begin{equation}
\nu_{\tilde{z}}(f_i)=\nu(f_i)\label{eq:nuztilde=nufi}
\end{equation}
for $i=1,\ldots,s$.

(2) Assume that $\init_\nu(\tilde{z}) \in \gr_\nu(K)$. 
Then there exists $\phi^* \in K$ such that, letting $z^*=\tilde{z}-\phi^*$, we have 
$\delta(f_i,z^*)=0$ for $i=1,\ldots,s$. 
\end{cor}

\noi\textit{Proof:} As before, (1) implies (2) by Remark \ref{rek:delta0}, so it is enough to prove (1). It follows from Lemma \ref{lem:ztildez} that if $z^* = z-\phi$ with $\nu(z^*) > \nu(z) = \nu(\phi)$, then 
\begin{equation}
\nu_{z^*}(g) \ge \nu_z(g)\quad\text{ for all }g \in K[z].\label{eq:monotonicity}
\end{equation}
We construct $\tilde z$ satisfying (\ref{eq:nuztilde=nufi}) recursively in $i$. Take the greatest integer $j\in\{0,\dots,s-1\}$ such that
\begin{equation}
\nu_z(f_i)=\nu(f_i)\quad\text{ for all }i\in\{1,\dots, j\}.\label{eq:nuz=nufi}
\end{equation}
Apply Lemma \ref{lem:gene} to $f_{j+1}$. We obtain an element $\tilde z=z-\phi$ such that  (\ref{eq:nuztilde=nufi}) holds with
$i=j+1$. Moreover, by (\ref{eq:nuz=nufi}) and (\ref{eq:monotonicity}), equality (\ref{eq:nuztilde=nufi}) also holds for all $i\le j$. This completes the proof by induction on $j$.
$\Box$
\medskip

\begin{deft}
	Let $S = \{v_1,\ldots,v_s\}$ be a finite subset of $ \Q\oplus (\Gamma\otimes_\Z \Q)$. The \textbf{convex hull} of $S$ is
	\[
	\left\lbrace t_1v_1 + \cdots + t_sv_s \in   \Q\oplus(\Gamma\otimes_\Z\Q)  \ \left|\ \sum_{i=1}^st_i =1, t_i \in \Q_{\ge0} \right.
	\right\rbrace.
	\]
\end{deft} 
Fix a polynomial $g = a_dz^d +a_{d-1}z^{d-1} \cdots+ a_0 \in K[z]$.
\begin{deft}\label{deft:NewtonP}
	The \textbf{Newton polygon} of $g$, denoted by $\Delta(g,z)$, is the convex hull of the set
	\[
	\bigcup_{i=0}^d \left\lbrace (i,\nu(a_i)) \right\rbrace \subset  \Q \oplus(\Gamma\otimes_\Z \Q).
	\]
\end{deft}

Let $g=a_d\prod\limits_{j=1}^dg_j$ be the factorization of $g$ into linear factors, with $g_j =z-\phi_j$.
\medskip

\textbf{From now till the end of the paper, assume that }$\nu(z)\ge0$.
\medskip

Let $L:=[(i,\epsilon),(j,\theta)] \subset\mathbb{Q} \oplus\Gamma$ be a segment with $i<j$. The \textbf{slope} of $L$, denoted $sl(L)$, is defined by \[sl(L):=\frac{\theta-\epsilon}{j-i}\in\Gamma\otimes_{\Z}\Q.\]
\begin{rek}
	The Newton polygon $\Delta(g_j,z)$  is the segment $[(1,0),(0,\nu(\phi_j))]$. 
	Now
\[
\Delta(g,z) = (0,\nu(a_d))+\sum_{j=1}^d \Delta(g_j,z)
\]
where the sum stands for the Minkowski sum. For each $j \in \{1,\ldots,d\}$, 
	$\Delta(g,z)$ has a side $L_j$ parallel to $[(1,0), (0,\nu(\phi_j))]$. 
\end{rek}

\begin{deft}
	In this situation, we say that $g_j$ \textbf{is attached to} $L_j$.
\end{deft}

\begin{rek}\label{privilegedornot}
	\begin{enumerate}
		\item Note that, for $\psi \in K$, we have $\nu(z) > \nu(\psi)$ if and only if $\nu(z) > \nu(z-\psi)$. 
		\item Take $g \in K[z]$ and let $h=z-\psi$ be a $\nu$-privileged factor of $g$. Then  $\nu(z) > \nu(\psi)$ if and only if
$\nu(z) > \nu(\phi_j)$ for all $j \in \{1,\ldots,d\}$. Indeed, ``if'' is trivial. ``Only if'' follows from (1) of this Remark:
\[
\nu(z) > \nu(\psi) \iff \nu(z) > \nu(z-\psi) \implies \nu(z) > \nu(z-\phi_j) \iff \nu(z) > \nu(\phi_j).
\]
	\end{enumerate}
\end{rek}

\begin{lem} \label{lem:deltanul}
	Take $g \in K[z]$ monic and let $h=z-\psi$ be a $\nu$-privileged factor of $g$.
	\begin{enumerate}				
		\item The following are equivalent:
		\begin{enumerate}
			\item $\nu(z) > \nu(\psi)$
			\item $\delta(g,z)=0$.
		\end{enumerate}	
		 \item If $\delta(g,z)=0$, every $\nu$-privileged factor
	 $h=z-\psi$ of $f$ is attached to the side of $\Delta(g,z)$ of the smallest slope. 
	\end{enumerate}
\end{lem}

\noi Proof: 1.(a)$\implies$(b) We number the $\phi_j$ in the increasing order of their values. Then, for each $j\in\{1,\dots,d\}$, we have
\[
\nu\left(\prod\limits_{q=1}^{d-j}\phi_q\right) \le\nu(a_j).
\]
Hence $\nu(z) > \nu(\psi) \Longleftrightarrow \nu(z) > \nu(z-\phi_j)$ for all
$j\in\{1,\ldots,d\}\Longrightarrow \nu(a_0) = \nu\left(\prod\limits_{i=1}^{d} \phi_i\right) < \nu\left(a_jz^j\right)$ for all
$j \in \{1,\ldots,d\} \Longleftrightarrow\delta(g,z)=0$.
\smallskip

\noi (b)$\implies$(a) We must show that the inequalities $\nu\left(\prod\limits_{i=1}^{d} \phi_i\right) < \nu\left(a_jz^j\right)$  for all  $j \in \{1,\ldots,d\}$ imply that
$\nu(z)  > \nu(\phi_i)$ for all $i \in \{1,\ldots,d\}$. 

We argue by contradiction. Number the $\phi_i$ so that 
\begin{equation} \label{eq:suitephi}
\nu(\phi_1) \le \cdots \le \nu(\phi_d).
\end{equation} Since $ \nu\left(\prod\limits_{i=1}^{d} \phi_i\right) < \nu(z^d)$, we have $\nu(\phi_1) < \nu(z)$. Assume that there exists
$i \in \{2,\ldots,d\}$ such that 
\begin{equation} \label{eq:compare}
\nu(\phi_i) \ge \nu(z)
\end{equation}
 and take the smallest such $i$. We have $a_{d-i+1} = \phi_1 \cdots \phi_{i-1} + \prod\limits_{(j_1,\ldots,j_{i-1})} \phi_{j_1} \cdots \phi_{j_{i-1}}$ where $(j_1, \ldots,j_{i-1})$ runs over all the $(i-1)$-tuples of elements of $\{1,\ldots,d\}$ other than $(1,\ldots,i-1)$. By (\ref{eq:suitephi}) and (\ref{eq:compare}), for each such $(i-1)$-tuple we have
\[
\nu(\phi_1 \cdots \phi_{i-1}) < \nu(\phi_{j_1} \cdots \phi_{j_{i-1}}).
\]
Hence $\nu(a_{d-i+1}) = \nu(\phi_1 \cdots \phi_{i-1})$ and 
\[
\nu\left(a_{d-i+1}z^{d-i+1}\right)= \nu\left(\phi_1 \cdots \phi_{i-1} z ^{d-i+1}\right) \le \nu\left(\prod_{\ell =1}^{d} \phi_\ell\right).
\]
This is a contradiction. This completes the proof of the first part of the Lemma.
\smallskip

2. By Part 1 of this Lemma, we have $\nu(z) > \nu(\psi)$. By Remark \ref{privilegedornot}, we have $\nu(z) > \nu(\phi_j)$, for all
$j\in\{1,\ldots,d\}$. Hence $\nu(z-\phi_j)<\nu(z)$ for all $j \in \{1,\ldots,d\}$. 

Now, $sl([1,0),(0,\nu(\phi_j)]) = -\nu(\phi_j)$. Thus $\{-\nu(\phi_1),\ldots,-\nu(\phi_d)\}$ is a complete list of slopes of sides of $\Delta(g,z)$ of the form $L_k$. The result follows immediately. $\Box$
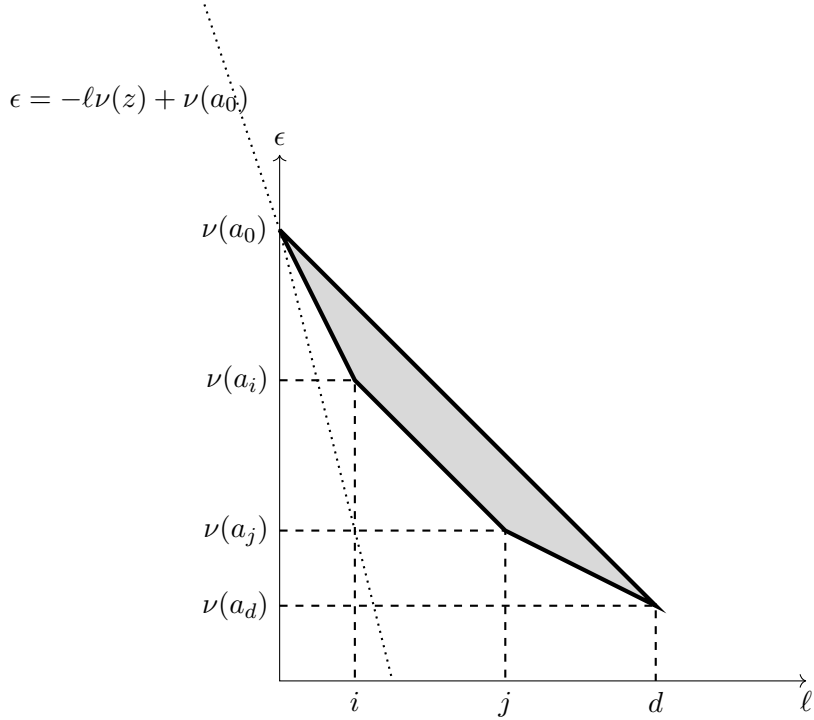
\medskip

\begin{tikzpicture}
\draw[<->] (7,0) node[below]{$\ell$} --(0,0)--(0,7) node[left,above]{$\epsilon$}; 
\draw[ultra thick,fill=LightGray] (0,6) node[left]{$\nu(a_0)$}--(1,4)--(3,2)--(5,1)--(0,6);
\draw[dashed,thick] (0,4) node[left]{$\nu(a_i)$}--(1,4)--(1,0) node[below]{$i$};
\draw[dashed,thick] (0,2) node[left]{$\nu(a_j)$}--(3,2)--(3,0) node[below]{$j$};
\draw[dashed,thick] (0,1) node[left]{$\nu(a_d)$}--(5,1)--(5,0) node[below]{$d$};
\draw[dotted,thick] (-1,9) -- (0,6) -- (1.5,0) node[above=220pt,left=50pt]{$\epsilon =  - \ell\nu(z) + \nu(a_0)$};
\end{tikzpicture}

\begin{rek}\label{rem:largeslope}
	We have $\delta(g,z) = 0$ if and only if $\delta(g_{red},z) = 0$. 
\end{rek}
\begin{prop} \label{prop:privbranch}
	Let $g=\prod\limits_{j=1}^d g_j$ be as above. Take an integer $k \in \{1,\ldots,d-1\}$. Let
$g^{(k)}=\prod\limits_{j=1}^{d-k}h_j$ be a 
	decomposition of $g^{(k)}$ into linear factors in $K[z]$ and let 
	$h$ be a $\nu$-privileged factor of $g^{(k)}$.  If 
	\begin{equation}
	\nu\left(g^{(i)}\right) \ge \nu(g)\quad \text{ for all }i \leq 
	k\label{eq:valeurderivee}
	\end{equation}
	then
	\begin{equation} \label{eq:inegval}
	\nu(h) > \nu(g_j)\ \qquad \text{ for } 1 \le j \le d.
	\end{equation}
\end{prop}

\noi\textit{Proof:} By Lemma \ref{lem:gene}, we can choose 
coordinates so that one of the following conditions holds:

(a) \begin{equation} \label{eq:deltanulle}
\delta(g,z)=0
\end{equation} or

(b) $\init_\nu z\notin gr_\nu K$.
\medskip

Next, we show that (b)$\implies$(a). Suppose (b) holds, that is, $\init_\nu z \notin \gr_\nu K$. By Remark \ref{rek:algebric}, we have
$\nu(g) = \nu_z(g)$.
\medskip

Let $\delta= \delta(g,z)$. We give a proof of (a) by contradiction. Assume that $\delta > 0$. Then
$(\delta-1,\nu(a_\delta))$ is a vertex of $\Delta(g',z)$. We have
\[
\nu(g) = \nu_z(g) = \min_{j \in \{0,\ldots,d\}}\left\{\nu\left(a_jz^j\right)\right\} = \nu\left(a_\delta z^\delta\right).
\]
On the other hand, 
\begin{align*}
	\nu\left(g'\right) &= \nu_z\left(g'\right) = \min_{j \in \{1,\ldots,d\}}(\nu(ja_jz^{j-1})) =  \min_{j\in \{1,\ldots,d\}}(\nu(a_jz^j)-\nu(z)) \\&= \min_{j \in \{1,\ldots,d\}}(\nu(a_jz^j))-\nu(z)) = \nu_z(g) -\nu(z) = \nu(g) -\nu(z)  < \nu(g).
\end{align*}
This contradicts (\ref{eq:valeurderivee}). This proves the implication (b)$\implies$(a). 
\medskip

Therefore it is sufficient to prove the Proposition under the assumption (a).  By Remark \ref{rek:delta0}, we have $\nu(g)=\nu_z(g)$. More precisely, we have
\begin{equation} \label{eq:nua0}
\nu(g) = \nu_z(g)= \nu\left(a_0\right) < \nu\left(a_jz^j\right)\quad\text{ for all } \ j\in\{1,\dots,d\}.
\end{equation}

\noi By the implication (b)$\implies$(a) of Lemma \ref{lem:deltanul}, formula (\ref{eq:nua0}) implies that
\[
\nu(z) > \nu(g_j).
\]
If $\nu(z) \le \nu(h)$, the 
proof is finished.
\smallskip

Hence we may assume that $\nu(z) > \nu(h)$. By Lemma \ref{lem:deltanul} (1),  we have
\[
\delta\left(g^{(k)},z\right) = 0.
\]
In particular,
\begin{equation} \label{eq:nudez}
\nu\left(g^{(k)}\right) =\nu(a_k).
\end{equation}

We have 
\begin{equation}
\nu(a_k)=\nu_z\left(g^{(k)}\right) =\nu\left(g^{(k)}\right) \ge \nu(g) = \nu_z(g)=\nu(a_0).\label{eqnuak>nua0}
\end{equation}
Let
\[
L = [ (0,\nu(a_0)),(\epsilon,\nu(a_\epsilon))]
\]
be the side of $\Delta(g,z)$ of minimal slope. Let
\[
\tilde{L}^{(k)} = [(0,\nu(a_k)),(\epsilon-k,\nu(a_\epsilon))]\subset\Delta\left(g^{(k)},z\right).
\]
Let $L^{(k)}$ be the side of $\Delta\left(g^{(k)},z\right)$ of minimal slope. We have $(0,\nu(a_k))\in\Delta\left(g^{(k)},z\right)$. 

Now let $g_j$ be a $\nu$-privileged factor of $g$. Using (\ref{eqnuak>nua0}) we obtain 
\begin{align*}
\nu(g_j)=-sl(L)=\frac{\nu(a_0) - \nu(a_\epsilon)}{\epsilon} < \frac{\nu(a_k)-\nu(a_\epsilon)}{\epsilon-k} \\
= -sl\left(\tilde{L}^{(k)}\right)\le-sl\left(L^{(k)}\right) = \nu(h).   
\end{align*}
This completes the proof of the Proposition. $\Box$
\medskip

\section{Semi-algebraicity of roots of polynomials}\label{semialgebraicity}

Fix a monic polynomial $g$ in one variable $z$ of degree $d$ whose coefficients are continuous semi-algebraic functions defined on a semi-agebraic set $U \subset R^n$.
\medskip

Recall the following definition from \cite{HI}:
\begin{deft} \label{deft:multival}
	A multiple-valued function $\mathcal{F}$ from a space $X$ to a space $Y$ will be called a \textbf{continuous }$n$\textbf{-valued function} from $X$ to $Y$ and will be denoted by $\mathcal{F}: X \to^n Y$ provided
	
	(i) to each $x \in X$, $\mathcal{F}$ assigns $m_x$ values $y_1,\ldots,y_{m_x}$ in $Y$ with associated multiplicities $k_i$ such that $\sum\limits_{i=1}^{m_x} k_i =n$;
	
	(ii) to each neighborhood $N(y_i)$ in $Y$ there corresponds a neighborhood $U(x)$ in $X$ such that for $z \in U(x)$ there are
$k_i$ values of $\mathcal{F}(z)$ in $N(y_i)$, counted with multiplicities.   
\end{deft}

\begin{thm} \label{thm:dvalue}
	The root $\mathcal{F}$ in $\sur{R}$ of $g(a)$, $a \in U$, is a continuous $d$-valued function from $U$ to $\sur{R}$.
\end{thm}

\begin{rek}
	This Theorem is well known in the case $\sur{R} = \C$, but the proofs found in the literature (\cite{Marden}, p. 3)) usually use complex analysis (namely, the characterization of the number of roots in a disk as a contour integral) and so are not applicable in our more general situation. The proof given below  uses only the Euclidean topology on $\sur{R}^d$ and some elementary facts about finite group actions on varieties.  
\end{rek}

\noi\textit{Proof of the Theorem:} Fix a point $a \in U$. Let $\psi(a)=(\psi_1(a), \ldots,\psi_d(a)) \in \sur{R}^d$ be the roots of $g(a)$, let $\ell \le d$ be the number of distinct roots among the $\psi_i(a)$. Renumbering the $\psi_i(a)$ if necessary, we may assume that there exist integers
\[
0=i_0 < i_1 < i_2 < \cdots <i_\ell = d
\]
such that: 

\noi (1) for all $j \in \{1,\ldots,\ell\}$ we have $\psi_{i_{j-1}+1}(a) = \psi_{i_{j-1}+2}(a) = \ldots = \psi_{i_j}(a)$

\noi (2) for all $j,j' \in \{1,\ldots,\ell\}$, $j\neq j'$, we have $\psi_{i_j}(a) \neq \psi_{i_{j'}}(a)$.
\smallskip

We want to check the definition of continuity at $a$. In other words, we must show that for any $\epsilon \in R$, $\epsilon > 0$, there exists $\delta \in R$ such that for all $a' \in U$, $|a'-a| < \delta$, there exists a numbering $\psi_1(a'),\ldots,\psi_d(a')$ of the roots of $g(a')$ such that $|\psi_j(a)-\psi_j(a')| < \epsilon$ for all $j \in \{1,\ldots,d\}$.

Let $S_d$ denote the symmetric group on $d$ elements. Let $S_d \times \sur{R}^d \to \sur{R}^d$ be the action of $S_d$ on
$\sur{R}^d$ by permutation of coordinates. Let $\sigma: \sur{R}^d \to \sur{R}^d/S_d \cong \sur{R}^d$ be the resulting quotient map, where the last isomorphism is obtained using the symmetric functions. Let
$b=(b_1,\ldots,b_d)$ be coordinates on the source $\sur{R}^d$. Then the symmetric functions $\sigma_1(b),\ldots,\sigma_d(b)$ form a natural coordinate system on the target $\sur{R}^d$. 
\begin{deft}
	Let $\xi=(\xi_1,\ldots,\xi_d)$ be a point of $\sur{R}^d$ and $\epsilon$ a strictly positive element of $\sur{R}$. The $\epsilon$-polydisk neighborhood of $\xi$ is the set
\[
D_\epsilon(\xi)=\left\{\left.\left(\xi'_1,\ldots,\xi'_d\right) \in \sur{R}^d\ \right|\ 
\left|\xi_j-\xi_{j'}\right| < \epsilon, j \in \{1,\ldots,d\}\right\}.
\]
\end{deft}
Continuity of the $d$-valued function $a \mapsto \psi(a)$ follows from the fact  that the map $\sigma$ is open at $\psi(a)$, that is, for each $\epsilon$-neighborhood  $D_\epsilon(\psi(a))$ of $\psi(a)$, its image $\sigma(D_\epsilon(\psi(a)))$ contains a
$\delta$-neighborhood $D_\delta(g(a))$ of $g(a)$. It remains to prove that $\sigma$ is open.

For $j \in \{1,\ldots,\ell\}$, put $k_j=i_j-i_{j-1}$. We consider the stabilizer of $\psi(a)$ in $S_d$: $Stab(\psi(a)) \cong S_{k_1} \times \cdots \times S_{k_\ell}$. Let \[c := \frac{d!}{k_1! \cdots k_\ell!}.\] Let $p_1,\ldots,p_c \in S_d$ be a set of representatives of the left cosets of $Stab(\psi(a))$. Then $\{p_1\psi(a),\ldots,p_c\psi(a)\}$ is the orbit of $\psi(a)$ under the action of $S_d$. 

Take $\epsilon$ sufficiently small so that all the $D_\epsilon(p_j \psi(a))$, $1 \le j\le c$, are disjoint. The group $S_d$ permutes the polydisks $D_\epsilon(p_j\psi(a))$, $1 \le j \le c$. Thus the action of $S_d$ on $\sur{R}^d$ restricts to an action on
$\coprod\limits_{j=1}^c D_\epsilon(p_j\psi(a))$. For each $j \in \{1,\ldots,c\}$, the action of $S_d$ on $\sur{R}^d$ induces an action of $p_jStab(p_j\psi(a))p_j^{-1}$ on $D_\epsilon(p_j\psi(a))$. We have the induced action of $S_d$ on
$\coprod\limits_{j=1}^c D_\epsilon(p_j\psi(a))/p_jStab(p_j\psi(a))p_j^{-1}$, having no fixed points. 

Consider the quotient $\sigma_{St}: D_\epsilon(\psi(a)) \to D_\epsilon(\psi(a))/Stab(\psi(a))$ of $D_\epsilon(\psi(a))$. In view of the above, it remains to show that $\sigma_{St}$ is open at $\psi(a)$. Write \begin{align*}
	D_\epsilon(\psi(a))= D_\epsilon(\psi_1(a),\ldots,\psi_{i_1}(a)) \times D_\epsilon(\psi_{i_1+1}(a),\ldots,\psi_{i_2}(a)) \times \cdots \\ \cdots \times D_\epsilon(\psi_{i_{\ell-1}+1}(a),\ldots,\psi_{i_\ell}(a)).
\end{align*}

Now, $S_{k_1} \times S_{k_2} \times \cdots \times S_{k_\ell}$ acts diagonally on $D_\epsilon(a)$, that is, given
$$
p=(p_1,\ldots,p_\ell) \in S_{k_1} \times S_{k_2} \times \cdots \times S_{k_\ell}
$$
and $(\lambda_1,\ldots,\lambda_\ell) \in D_\epsilon(\psi(a))$, we have
$p\lambda=(p_1\lambda_1,p_2\lambda_2,\ldots,p_\ell \lambda_\ell)$. Thus to show that $\sigma_{St}$ is open at $\psi(a)$, it is sufficient to show that, for each  $j \in \{1,\ldots,\ell\}$, the quotient map of $D_\epsilon(\psi_{i_{j-1}+1}(a),\ldots,\psi_{i_j}(a))$ by $S_{k_j}$ is open at $(\psi_{i_{j-1}+1}(a),\ldots,\psi_{i_j}(a))$.

This reduces the proof of the Theorem to the following lemma (where we have replaced $k_j$ by $d$ and $(\psi_{i_{j-1}+1}(a),\ldots,\psi_{i_j}(a))$ by the origin in $\sur{R}^d$). 
$\Box$  

\begin{lem} \label{lem:origin}
	The map $\sigma$ is open at the origin of $\sur{R}^d$. 
\end{lem}	
	
\noi Proof: Take a strictly positive $\epsilon \in R$. Let $\delta = \min\left\{\left(\frac{\epsilon}{d!}\right)^{d!},\frac12\right\}$. Take a point $a=(a_1,\ldots,a_d) \in D_\delta(0)$ and a point $b=(b_1,\ldots,b_d) \in \sigma^{-1}(a)$. We want to prove that \begin{equation} \label{eq:depsilon} b \in D_\epsilon(0).
\end{equation}
We will prove (\ref{eq:depsilon}) by induction on $d$. For $d=1$ the result is obvious. Assume that $d>1$ and (\ref{eq:depsilon}) holds with $d$ replaced by $d-1$. 

Without loss of generality, assume that $|b_1| \le |b_2| \le \cdots \le |b_d|$. Since
$$
|a_d| = |\sigma_d(b)| =\left|\prod\limits_{j=1}^d b_j \right| < \delta,
$$
we have $|b_1| < \delta^{1/d} < \epsilon$. 

Let 
\[
z^{d-1} + c_1 z^{d-2} + \cdots + c_{d-1} = \frac{z^d+a_1z^{d-1} + \cdots + a_d}{z-b_1}.
\]
Each $c_i$ is a sum of at most $d$ terms, each of which has absolute value strictly less than $\delta^{1/d}$. Thus
\[
|c_i| < d\delta^{1/d} \le d\left(\frac{\epsilon}{d!}\right)^{(d-1)!} \le \left(\frac{\epsilon}{(d-1)!}\right)^{(d-1)!}.
\]
By the induction assumption, we have $|b_i| < \epsilon$ for $i \in \{1,\ldots,d\}$ as desired. 

This completes the proof of Lemma \ref{lem:origin} and Theorem \ref{thm:dvalue}.
$\Box$
\medskip

Assume that for $a \in U$ the number of real roots of $g(a)$, counted with multiplicities, is independent of $a$. Let $s$ denote this number of real roots, common for all $a$. We define the following functions $\phi_i:U \to R$, $i \in \{1,\ldots,s\}$. For each $a \in U$, define $(\phi_1(a),\ldots,\phi_s(a))$ to be the $s$-tuple of real roots of $g(a)$ in $R$ arranged in the increasing order.

\begin{cor} \label{cor:continu}
	The functions $\phi_i$ are continuous.
\end{cor}

\noi First of all, we recall two Lemmas from \cite{HI}, generalizing the first one from $\R$ to any closed real field $R$:

\begin{lem} (Lemma 1, p. 431) 
	For any $\mathcal{F} : X \to^n R$, the least value $f(x)$ of $\mathcal{F}(x)$ is a continuous function. 
\end{lem}

The proof is exactly the same as the one given in  \cite{HI}, but we prefer to reproduce it here for the sake of completeness.  
\medskip

\noi Proof : Let $x_0 \in X$. By the property (ii) of Definition \ref{deft:multival}, for any $\epsilon \in R$, $\epsilon  >0$, there exists a neighborhood $U$ of $x_0$ such that for all $z \in U$, all $n$ values of $\mathcal{F}(z)$ are greater than $f(x_0)-\epsilon$. Similarly, using the property (ii) again, there exists a neighborhood $V$ of $x_0$ in which at least one value of $\mathcal{F}(z)$ is less than $f(x_0)+\epsilon$. Hence, if $z \in U\cap V$, we have $f(x_0)-\epsilon < f(z) < f(x_0)+ \epsilon$, so $f$ is continuous at $x_0$.  $\Box$

\begin{lem} (Lemma 3, p. 432)
	If $\mathcal{F}:X \to^n Y$  has always exactly $m$ values in an open or closed subspace $W$ of $Y$, then the restriction of the values of $\mathcal{F}(x)$ to $W$ defines a continuous multi-valued function $\mathcal{F}':X \to^m W$.
\end{lem}
 \medskip

\noi\textit{Proof of Corollary \ref{cor:continu}:} We proceed by induction on $d$. For $d=1$ the Corollary is obvious. Together with Theorem \ref{thm:dvalue}, the preceding Lemmas imply that $\phi_1$ is continuous. Now, $z-\phi_1$ divides $g$, so
$\tilde{g}=\frac{g}{z-\phi_1}$ is a polynomial of degree $ d-1$, whose coefficients are continuous functions. We may apply the induction hypothesis to $\tilde{g}$ and conclude that all the real roots of $\tilde{g}$ are continuous. This completes the proof of the Corollary.
$\Box$
\medskip

For $1\le i\le s$ let $G_i:=graph(\phi_i)=\{(a,\phi_i(a))\ |\ a\in U\}$.
\medskip

The next Proposition is a special case of a general result by H. Delfs and M. Knebusch in 1981 : Lemma 1.1 of \cite{DK} which proves the semi-algebraicity of the roots of a polynomial whose coefficients are semi-algebraic functions of $x$.

\begin{prop}
	Let the notation be as above. Assume, in addition, that the coefficients of $g$ are semi-algebraic functions of $x$. Then for each $i \in \{1,\ldots,s\}$ the set $G_i$ is semi-algebraic. The functions $\phi_i$ are semi-algebraic (in other words, for each $i$ the linear polynomial $z-\phi_i$ is a real branch of $g$ over $U$). The projection $\pi|_{G_i}: G_i \to U$ is a semi-algebraic homeomorphism. 
\end{prop}

\section{Value of a semi-algebraic function at a point of the real spectrum }\label{Value}

Recall the following notation from the Introduction. We consider the Euclidean space $R^{n+1}$ with coordinates $(x,z)$, where $z$ is a single variable and $x=(x_1,\dots,x_n)$. Let $\pi: R^{n+1} \to R^n$ be the projection onto the $x$-space and 
\[
\tilde{\pi} : \sper\ B[z] \to \sper\ B
\]
be the corresponding morphism of real spectra. Let  $D\subset R^{n}$ be a connected semi-algebraic set; consider the cylinder
$C=\pi^{-1}(D)$.

Let $\gamma \in \sper\ B[z]$ be a point of the real spectrum of $B[z]$, $\gamma \in \tilde{C}$ (where $\tilde C$ is the cylinder in $\sper\ B[z]$, corresponding to $C$); let  $\gamma_0 := \tilde{\pi}(\gamma) \in \tilde{D}$. Recall that for $g \in B[z]$, 
$g(\gamma)$ denotes the natural image of $g$ in $B[z](\gamma)$. Let $h=z-\phi$ be a real branch of $g$ over $D$ (cf. Definition \ref{def:realbranch}). 

Following \cite{BCR} \S 7.3, we can associate to $\phi$   a map  $\tilde{\phi}: \tilde{D} \to \coprod\limits_{\eta_0 \in \tilde{D}} \sur{B(\eta_0)}_r$. Let $h(\gamma)$ be the image of $z-\tilde{\phi}(\gamma_0)$ under the natural homomorphism $\sur{B(\gamma_0)}_r[z] \to\sur{B[z](\gamma)}_r$.
\medskip

Let $U$ be a connected semi-algebraic set such that the number of real roots of $g$ over $U$, counted with or without multiplicity, is constant. Let $s$ be the number of distinct real roots of $g$. 
\medskip

We have the functions $\tilde{\phi}_i : \tilde{U} \to \coprod\limits_{\eta_0 \in \tilde{U}} \sur{B(\eta_0)}_r$ and continuous semi-algebraic functions $\phi_i:U \to R$, $i \in \{1,\ldots,s\}$ such that, for each closed point $a \in U \subset \tilde{U}$, we have $\phi_i(a) = \tilde{\phi}_i(a)$ (here we identify $a \in U$ with its natural image under the natural injection $U \subset \tilde{U}$ and $\sur{B(a)}_r$ with $R$).
\medskip

 For each $\eta_0 \in \tilde{U}$ (resp. $a \in U$), $(\tilde{\phi}_1(\eta_0), \ldots,\tilde{\phi}_s(\eta_0))$ (resp. $(\phi_1(a),\ldots,\phi_s(a))$) is the $s$-tuple of distinct real roots of $g(\eta_0)$ in $\sur{B(\eta_0)}_r$ (resp. real roots of $g(a)$ in $R$) arranged in the increasing order. 
\medskip

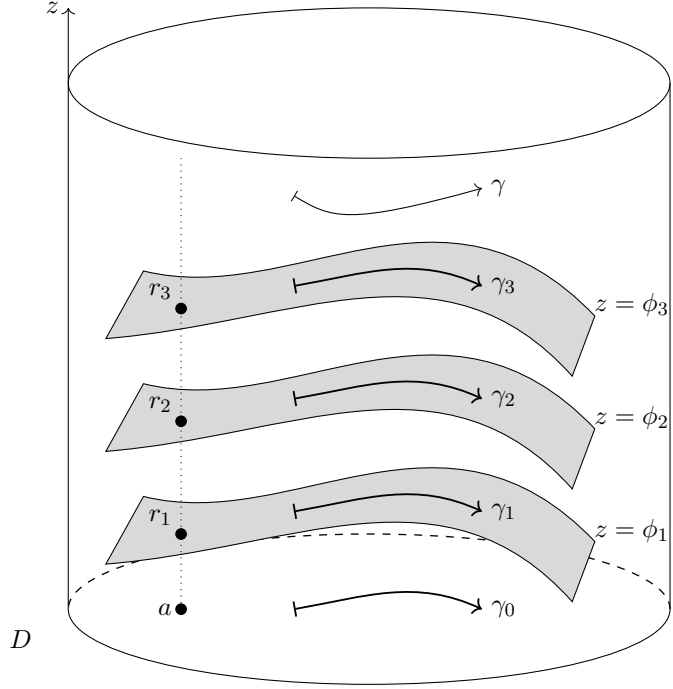
\begin{center}
  \tikzset{
    partial ellipse/.style args={#1:#2:#3}{
      insert path={+(#1:#3) arc (#1:#2:#3)}
    }
  }
  \begin{tikzpicture}[every node/.style={scale=0.9}]
    \draw[->] (8,0) -- (8,7) (0,0) -- (0,8) node[left] {$z$};
    \draw (4,0) [partial ellipse=180:360:4 and 1];
    \draw[dashed,line width=.5pt] (4,0) [partial ellipse=0:180:4 and 1] node[left=20pt,below=5pt]{$D$};
    \draw  (4,7) ellipse (4 and 1) ;

    \foreach \n in {1,2,3}{  
      \draw[line width=.4pt,fill=LightGray,shift={(0,1.5*(\n-1))}]
        (0.5,0.6) -- (1,1.5) .. controls (3,1) and
        (5,3) .. (7,0.9) node[above right=-4pt]{$z=\phi_\n$}
        -- (6.7,0.1) .. controls (5,2) and (3,0.8) .. (0.5,0.6) -- cycle;
      \draw[|->,line width=.7pt,shift={(0,1.5*(\n-1))}] 
        (3,1.3) to[out=10,in=155] (5.5,1.3) node[right]{$\gamma_\n$};
      \draw[fill=black,shift={(0,1.5*(\n-1))}] 
        (1.5,1) circle (2pt) node[above left] {$r_\n$};
    }

    \draw[|->,line width=.7pt,shift={(0,-1.3)}] 
    (3,1.3) to[out=10,in=155] (5.5,1.3) node[right]{$\gamma_0$};

    \draw[|->] (3,5.5) .. controls (3.5,5.2) and (3.7,5.1) .. (5.5,5.6)
    node[right]{$\gamma$};

    \draw[fill=black] (1.5,0) circle (2pt) node[left]{$a$} ;
    \draw[dotted] (1.5,0) -- (1.5,6);    
  \end{tikzpicture}
\end{center}

For $g \in B$, let $g(\gamma_0) \in B(\gamma_0)[z]$ be the polynomial obtained from $g$ by replacing all the coefficients by their images in $B(\gamma_0)$. This way we have associated to each $\gamma_0 \in \tilde{D}$ a collection of $s$ real factors
$z-\tilde{\phi}_1(\gamma_0),\ldots,z-\tilde{\phi}_s(\gamma_0)$ of $g(\gamma_0)$. Conversely, given a real factor
$z-\tilde{\phi}(\gamma_0)$ of $g(\gamma_0)$, there exists $j \in \{1,\ldots,s\}$ such that $\tilde{\phi}(\gamma_0) = \tilde{\phi}_j(\gamma_0)$. 
\medskip

\begin{rek} \label{rek:bijection}
	The above results show that, for any $\gamma_0 \in \tilde{D}$, there is a natural order-preserving bijection between real branches of $g$ over $D$ and real $\gamma_0$-branches.
\end{rek}

\section{Real and imaginary parts of branches}\label{Sec:reel}

 Denote by $A$ the ring $B[z]=R[x,z]$. Let $\gamma \in \sper\ A$. The point $\gamma$ determines morphisms
\[
A[\gamma]= \frac{A}{\mathfrak{p}_\gamma} \hookrightarrow \sur{A(\gamma)}_r\hookrightarrow \sur{A(\gamma)}.
\]
Now, $\sper\ \sur{A(\gamma)}_r$ consists of a single point $\sur{\gamma}$. The valuation associated to this point is a natural extension of $\nu_\gamma$ to $\sur{A(\gamma)}_r$ which we denote by $\nu_{\sur{\gamma}}$.

We can view $\sur{A(\gamma)}$ as the field extension $\sur{A(\gamma)}_r(i)$ where $i^2=-1$. Thus any $\xi \in \sur{A(\gamma)}$ can be written as $\xi = u+iv$ where $u = Re\; \xi,v = Im\; \xi \in \sur{A(\gamma)}_r$. 

The purpose of this section is to study the extension of $\nu_{\bar\gamma}$ to $\sur{A(\gamma)}$; by abuse of notation this extension will also be denoted by $\nu_{\bar\gamma}$. The main result is
\begin{prop} \label{lem:reim}
The valuation $\nu_{\sur{\gamma}}$ admits a unique extension to $\sur{A(\gamma)}$, also denoted by $\nu_{\sur{\gamma}}$. This extension is characterized by the fact that for each $h \in \sur{A(\gamma)}$ we have 
\begin{equation} \label{eq:nu-alpha}
 \nu_{\sur{\gamma}}(h) = \min\{\nu_{\sur{\gamma}}(Re\, h), 
\nu_{\sur{\gamma}}(Im\, h) \}.
\end{equation}
\end{prop}

\noi\textit{Proof:} Let $\nu'$ be an extension of $\nu_{\sur{\gamma}}$  to $\sur{A(\gamma)}$. Take an 
element $h \in \sur{A(\gamma)}^*$. Since $\nu'(i)=0$, we have 
\begin{equation} \label{eq:extens}
 \nu'(h) \ge\min\{\nu_{\sur{\gamma}}(Re\, h), \nu_{\sur{\gamma}}(Im\, h)\}.
\end{equation}
Similarly, letting $\bar h$ denote the complex conjugate of $h$, we have
\begin{equation} \label{eq:extens2}
 \nu'(\sur{h}) \ge \min\{\nu_{\sur{\gamma}}(Re\, h),
\nu_{\sur{\gamma}}(Im\, h)\}.
\end{equation}
To prove that (\ref{eq:extens}) and (\ref{eq:extens2}) are equalities, we may assume that $Re\, h 
\neq 0$ and $Im\, h \neq 0$.  Since $(Re\, h)^2 >_{\sur{\gamma}} 0$ and $(Im\, h)^2 
>_{\sur{\gamma}} 0$, we have
\[
(Re\, h)^2 + (Im\, h)^2 >_{\sur{\gamma}} (Re\, h)^2.
\]
and
\[
(Re\, h)^2 + (Im\, h)^2 >_{\sur{\gamma}} (Im\, h)^2.
\]
Hence 
\begin{equation}\label{eq:extens3}
 \nu_{\sur{\gamma}}\left((Re\, h)^2 + (Im\, h)^2\right) \le \nu_{\sur{\gamma}}\left((Re\, h)^2\right)
\end{equation}
and
\begin{equation}\label{eq:extens4}
 \nu_{\sur{\gamma}}\left((Re\, h)^2 + (Im\, h)^2\right) \le \nu_{\sur{\gamma}}\left((Im\, h)^2\right),
\end{equation}
so that
\begin{equation}\label{eq:extens5}
 \nu_{\sur{\gamma}}\left((Re\, h)^2 + (Im\, h)^2\right) \le
\min\left\{\nu_{\sur{\gamma}}\left((Re\, h)^2\right), \nu_{\sur{\gamma}}\left((Im\, h)^2\right)\right\}.
\end{equation}
Combining (\ref{eq:extens}), (\ref{eq:extens2}) and (\ref{eq:extens5}), we obtain
\begin{align*}      
\min\left\{\nu_{\sur{\gamma}}\left((Re\, h)^2\right), 
\nu_{\sur{\gamma}}\left((Im\, h)^2\right)\right\} & = 2 \min\left\{\nu_{\sur{\gamma}}\left(Re\, 
h\right),\nu_{\sur{\gamma}}\left(Im\, h\right)\right\} \le \nu'(h) + \nu'(\sur{h}) \\ & = 
\nu'(h\sur{h}) = \nu_{\sur{\gamma}}((Re\, h)^2 + (Im\, h)^2) \\ &\le 2 
\min\{\nu_{\sur{\gamma}}((Re\, h),\nu_{\sur{\gamma}}(Im\, h)\}.
\end{align*}

\noi Thus all the inequalities are equalities. Hence \[\nu'(h) = \nu'(\sur{h}) = 
\min\{\nu_{\sur{\gamma}}(Re\, 
h),\nu_{\sur{\gamma}}(Im\, h)\}.\] 
$\Box$ \medskip

\begin{cor} \label{cor:no1}
Let $h \in \sur{A(\gamma)}$. We have $\nu_{\sur\gamma}(Im\, h) \ge 
\nu_{\sur\gamma}(h)$. 
\end{cor}

\begin{cor} \label{cor:nuIm}
	Let $f \in A$. Let $h_\gamma$ be a $\gamma$-privileged 
	branch of $g$. For every $\gamma_0$-branch $\tilde{g}$ of $g$ we have 
	$\nu_{\sur\gamma}(Im\, h_\gamma) \ge \nu_\gamma(\tilde{g})$. 
\end{cor}

\noi\textit{Proof:} This is an immediate consequence of Corollary \ref{cor:nuIm} and the fact that $h_\gamma$ is 
$\gamma$-privileged. $\Box$

\begin{cor}
	Assume that $\nu_\gamma\left(g^{(k)}\right) \ge 
	\nu_\gamma(g)$. Let $h_\gamma$ be an $\gamma$-privileged branch of 
	$g^{(k)}$ and $\tilde{g}$ a $\gamma_0$-branch of $g$. Then $\nu_\gamma(Im\, h_\gamma) \ge 
	\nu_\gamma(h_\gamma) > \nu_\gamma(\tilde{g})$. 
\end{cor}

\noi\textit{Proof:} This follows from Corollary \ref{cor:nuIm} and Proposition \ref{prop:privbranch}. $\Box$
\medskip

\begin{rek} For any branch $h : D  \to \sur{R}$, we can define two continuous functions $Re\; h : D \to R$ (resp. $Im\; h : D \to R$) by setting, for each $a \in D$,
\[
(Re\; h)(a) = Re(h(a))
\]
(resp. $(Im\; h)(a) = Im(h(a))$). There exists a connected semi-algebraic set $U \subset D$ containing $\gamma_0$ such that $Re\; h$  and $Im\; h$ are real branches over $U$. We can apply the preceding construction to define, for any $\gamma \in \tilde{C}$, $Re\; h(\gamma)$ and $Im\; h(\gamma)$ as elements of $\sur{A(\gamma)}_r$. It is precisely in this situation that the results of this section will be most frequently applied.
\end{rek}

\section{Rolle-type theorems for the real spectrum}\label{Sec:Rolle}

In order to prove Theorem  \ref{thm:main}, we start with some lemmas.

Assume that $(0,0) \in \sur{D} \setminus D$.
\medskip

 Let $h_1 = z-\phi_1$ and $h_2=z-\phi_2$ be real branches over $D$ for some 
 connected semi-algebraic set $D \subset R^n$, 
 $\phi_1(0)=\phi_2(0)=0$. Let $\gamma$ be a point of the real spectrum 
 $\sper\ B[z]$.

\begin{deft}
	If $0 < h_1(\gamma) < h_2(\gamma)$ or 	$h_2(\gamma) < h_1(\gamma) <0$, we say that $h_1$ \textbf{lies between} $\gamma$ and $h_2$. 
\end{deft}

\begin{lem} \label{lem:conv2}
	If $h_1$ lies between $\gamma$ and $h_2$, then $\nu_\gamma(h_1(\gamma)) \geq 
\nu_\gamma(h_2(\gamma))$. 
\end{lem}

\noi\textit{Note:} geometrically  $0<h_1(\gamma) \leq h_2(\gamma)$ means 
that the curvette $\gamma$ is on the same side of the 
hypersurfaces $h_1=0$ and $h_2=0$. 
\medskip
 
\noi\textit{Proof:} Recall that by definition of $\nu_\gamma$ the valuation ring $R_\gamma$ is given by
\[
R_\gamma = \left\{\left.x \in \sur{B[z](\gamma)}_r \ \right|\ \exists y \in B[z][\gamma]\text{ such that }\  |x| \le |y|\right\}.
\]
By hypothesis, $|h_1| < |h_2|$, so $\left|\frac{h_1}{h_2}\right| <1$. 
Hence $\frac{h_1}{h_2} \in R_\gamma$, from which we deduce that $\nu_\gamma(h_1) \ge 
\nu_\gamma(h_2)$.
 $\Box$
\medskip

\begin{deft}
	The ring of Puiseux series $R[[t^\Q]]_{\text{Puisseux}}$ is the ring of generalized series over $R$ whose exponents are non-negative rational numbers with bounded denominators. Denote its quotient field by  $R((t^\Q))_{\text{Puisseux}}$
\end{deft}

\textbf{We will now restrict attention to points $\gamma \in C$ defined by Puiseux 
series}. In the next section we will show how to reduce the general case 
to the case of such semi-curvettes. 
\medskip

For $g \in A$ and $\gamma \in \sper\ A$ given 
by a curvette $\gamma(t) = (x_1(t),\ldots,$ $x_n(t),z(t))$, where $x_j(t)$ and $z(t)$ are generalized power series over $R$ whose exponents are rational with bounded denominators, put $\sur{g}(z,t)=g(z,x(t))$. This is a polynomial in $z$ over the ring of Puiseux series in $t$. As before, let $\gamma_0=\pi(\gamma)$. Fixing a curvette as above is equivalent to fixing a homomorphism
 $L:A[\gamma] \to R[[t^\Q]]_{\text{Puiseux}}$.  
 Below we will be interested in the restriction of $L$ to $B[\gamma_0]$, which induces a homomorphism $L_0 : B[\gamma_0][z] \to R[[t^\Q]]_{\text{Puiseux}}[z]$. Since the field $R((t^\Q))_{\text{Puiseux}}$ is real closed, the homomorphisms $L$ and $L_0$ extend naturally to homomorphisms $\iota
\medskip: \sur{A(\gamma)}_r \to R((t^\Q))_{\text{Puiseux}}$ and $\iota_0 : \sur{B(\gamma_0)}_r[z] \to R((t^\Q))_{\text{Puiseux}}$ respectively.

Now let $g=z^d+ a_{d-1}(x)z^{d-1} + \cdots + a_1(x) z+a_0(x)$ be a monic polynomial in $A=B[z]$, where $a_i(x)$, $i \in  \{0,\ldots,d-1\}$, are elements of $B$. Let $h_1,h_2$ be two 
real branches of $g$ over $D$. Write $h_i=z-\phi_i$, $i=1,2$ such that $\phi_1(0)=\phi_2(0)=0$.  
\medskip

For each $a \in \N$, the map $t \mapsto t^a$ 
induces a homeomorphism of the half-plane $\{t \ge 0\}$ onto itself. 
We can choose a positive integer $a$ such that, replacing 
$t$ by $t^a$ in 
$\sur{g}$, we may assume that 
$\sur{g}$ is a formal power series over $R$ with integer exponents. By the results of \S\ref{Value}, we can canonically associate to $\phi_i$ an element $\tilde{\phi}_i \in \sur{B(\gamma_0)}_r$. Let $\sur{\phi}_i= \iota_0(\tilde{\phi}_i)$. 
Up to performing a new change $t \mapsto t^a$, we may assume that all $\sur{g}$ 
and $\sur{\phi_i}$ are series with integer exponents.
\smallskip

Write $\sur{h}_1 = z-\sur{\phi}_1(t) = z -\sum\limits_{i=1}^\infty b_it^i$ and $\sur{h}_2 = z-\sur{\phi}_2(t) = z-\sum\limits_{i=1}^\infty c_it^i$. Let $s= \max\{i \in \N\ |\ b_i=c_i\}$. 

Replacing $z$ by $z-\sum\limits_{i=1}^{s+1} c_it^i$, we may assume that $c_{s+1}=0$ and 
$c_i=b_i=0$ for $i \leq s$. Also, up to interchanging $\phi_1$ and $\phi_2$, we may assume that $b_{s+1}> 0$. Let
$\kappa = \min\{j \in \N \ |\ c_{s+j} \neq 0\}$.
\medskip

For $h=\sum\limits_{i=0}^q e_iz^i \in R[[t]][z]$, let 
\begin{align} \label{eq:defval}
 \nu_s(h) &= \min\{(s+1)i + \nu_\gamma(e_i)\ |\ e_i \neq 0\} \\
\init_s h &= \sum_{(s+1)i+\nu_\gamma(e_i)=\nu_s(h)} e_iz^i.                
\end{align}

In order not to overload the notations, from now to the end of the section, we will write $g$ instead of $\sur{g}$.
 
\begin{lem} \label{lem:equid2}\textbf{(equidistance)}  (1) Let $h \in R[[t]][z]$ 
be a polynomial such that
\[
\init_s h = \init_s g'.
\]
There is at least one real factor of $h$ of the form $z-d_{s+1}t^{s+1} + h.o.t.$ with
\[
0<d_{s+1} <b_{s+1}.
\] 

(2) For each such factor $v$, there is no unique minimum 
among the three values $\nu_\gamma(h_1)$, $\nu_\gamma(h_2)$, $\nu_\gamma(v)$. 

Equivalently, we  have 
\begin{align*}
 \min\{\nu_\gamma(h_2),\nu_\gamma(v) \}  \le  \nu_\gamma(h_1),&\
\min\{\nu_\gamma(h_1),\nu_\gamma(v) \} \le  \nu_\gamma(h_2),\\
\min\{\nu_\gamma(h_1),\nu_\gamma(h_2) \}  & \le \nu_\gamma(v).
\end{align*}

\end{lem}

\noi\textit{Proof:} Let  $g=\sum\limits_{j,k}c_{jk}z^jt^k \in R[[t]][z]$. We view $R((t))$ as a valued field with the $t$-adic valuation. Let \[\Delta(g) 
= 
\text{ convex hull } \left\{  \bigcup_{c_{jk} \neq 0}(\{(j,k)\} ) 
\right\}.\] be the Newton polygon of $g$ (see Definition \ref{deft:NewtonP}).  

Let \begin{equation} \label{eqn:ell}
L=[(j_0,k_0),(j_1,k_1)], \text{ with }j_0<j_1,\ k_0> k_1
\end{equation}, be an edge of 
$\Delta(g)$ with strictly negative slope. The initial form of $g$ with 
respect to $L$ is defined to be
\[
\text{in}_L(g)  = \sum_{(j,k) \in L}c_{jk} z^jt^k.
\]
Let $L'=[(j_0-1,k_0),(j_1-1,k_1)]$ if $j_0>0$ and 
$L'=\left[(0,k_0-\frac{k_0-k_1}{j_1}),(j_1-1,k_1)\right]$ otherwise. Note that, if $j_1 \geq 
2$ (resp. $j_1=1$), $L'$ is a side (resp. a vertex) of $\Delta(g')$ and 
$(\text{in}_L(g))'=\text{in}_{L'}(g')$. \smallskip

By construction, $\Delta(g)$ has an edge $L$ with slope $-(s+1)$. Let the notations be as in (\ref{eqn:ell}). Let 
$\tilde{g}$ denote the polynomial in one variable $u$ such that 
\[\tilde{g}(\frac{z}{t^{s+1}}) = 
\frac{\text{in}_L(g(z,t))}{t^{(s+1)j_0+k_0}}.\]   

Since $\{g=0\}$ has a factor of the form 
$z-c_{s+\kappa}t^{s+\kappa} + \cdots$ with $s+\kappa > s+1$, $L$ cannot be the 
leftmost edge of $\Delta(g)$. In particular, $j_0 \ge 1$, so $\tilde{g}(0)=0$. 
As $\tilde{g}(b_{s+1})=0$, by Rolle's theorem (see for instance \cite{BCP}), there is at least one root $d_{s+1}$ of 
$\tilde{g}'$ in the open interval $(0,b_{s+1})$. This proves (1). 

Now fix one such $d_{s+1}$ and let $v$ be a factor of $\{h=0\}$ of the form 
$v=z-d_{s+1}t^{s+1} + h.o.t$.  It remains to prove that $v$ satisfies the 
conclusion of (2) of the Lemma.
 
Consider the set $\{\nu_\gamma(z-ct^{s+1}) \ |\ c \in R \}$. We have 
$\#\{\nu_\gamma(z-ct^{s+1}) \ |\ c \in R \} \le 2$. In other words, either $ 
\nu_\gamma(z-ct^{s+1})$ is constant for all $c \in R$ or there exists a unique $c^* \in R$ such that $\nu_\gamma(z-c^*t^{s+1}) >\nu_\gamma(z-ct^{s+1})$ for all $c \in 
R \setminus \{c^*\}$. 

If $\#\{\nu_\gamma(z-ct^{s+1}) \ |\ c \in R \} = 1$ or $c^* \notin 
\{0,b_{s+1},d_{s+1}\}$, then 
$\nu_\gamma(h_1)=\nu_\gamma(h_2)=\nu_\gamma(v) 
= \min \{\nu_\gamma(z),(s+1)\nu_\gamma(t)\}$ and (2) of the Lemma holds.

If $\#\{\nu_\gamma(z-ct^{s+1}) \ |\ c \in R \} = 2$ and $c^*$ coincides with 
one of $0,b_{s+1},d_{s+1}$, then the corresponding factor among 
$h_1,h_2,v$ has strictly greater $\nu_\gamma$-value than the other two, whose 
$\nu_\gamma$-values are equal. This completes the proof of the Lemma.  
$\Box$
\medskip

\begin{deft}
 Let $\alpha, \beta \in \sper\ A$, let $h$ be a real branch over $D$. 
 We say that $h$ is between $\alpha$ and $\beta$ if $h(\gamma) >0$ and $h(\beta) <0$ or vice versa. 
\end{deft}

\begin{deft} Let $z-\psi, z-\phi_1$ be real branches over $D$ and $z-\phi_2$ a $\gamma_0$-branch, not necessarily real. We say that $\psi$ lies between $\phi_1$ and 
$\phi_2$ if 
\begin{equation}
\phi_1(\gamma_0) \le \psi(\gamma_0) \le Re(\phi_2) \text{ or } Re(\phi_2) \le \psi(\gamma_0) \le \phi_1(\gamma_0).
\end{equation}
\end{deft}

\begin{lem} (generalized Rolle's Theorem) \label{lem:equid3}
Let $\phi_1 \in \sur{R(x)(\gamma_0)}_r$ and $\phi_2 \in \sur{R(x)(\gamma_0)}$ 
be roots of $g^{(k)}(\gamma_0)$ such that $\nu_\gamma(h_1) < 
\nu_\gamma(h_2)$ where $h_i$ stands for $z-\phi_i$. 

There exists a real root  $v \in \sur{R(x)(\gamma_0)}_r$ of 
$g^{(k+1)}(\gamma_0)$ between $\phi_1$ and $\phi_2$ such that \[\nu_\gamma(z-v) 
=\nu_\gamma(h_1) < \nu_\gamma(h_2).\]
\end{lem}

\noi Proof : Let \[\tilde{h}=g^{(k)}\frac{Re(h_2)^2}{h_2 \sur{h_2}}.\]
\smallskip

Note that \begin{equation} \label{eqn:nu}
           \nu_\gamma(h_1) < \nu_\gamma(h_2)
          \end{equation}
 implies that, after a 
change of variables, we may assume that $h_1=z-c_{s+1}t^{s+1} + h.o.t.$ and 
$h_2= z-b_{s+\kappa}t^{s+\kappa} + h.o.t.$ where $\kappa > 1$. Another 
consequence of (\ref{eqn:nu}) is : \begin{equation} 
\label{eq:inegval2} \nu_\gamma(h_1) = 
(s+1)\nu_\gamma(t) < \min \{\nu_\gamma(z),(s+\kappa)\nu_\gamma(t) 
\}\end{equation}
                                                                 
By Corollary \ref{cor:no1} and (\ref{eq:inegval2}), we have
$z=\init_s(h_2) = \init_s(Re(h_2))$ and the same for $\sur{h_2}$. Hence 
$\init_s(g^{(k)})=\init_s(\tilde{h})$. Hence 
\begin{equation} \label{eq:init}
 \init_s(g^{(k+1)}) = \init_s(\tilde{h}').
\end{equation}

Let $L$ be the side of the Newton polygon of $\tilde{h}$ whose slope is 
$-(s+1)$. Let $(j_0,k_0)$ be an integer point on $L$.  

Consider the polynomial $\tilde{g}(u)$ in one variable, defined by 
\[\tilde{g}(\frac{z}{t^{s+1}}) = \frac{\init_{L} 
\tilde{h}}{t^{(s+1)j_0+k_0}}.\]

Let $L'$ be the side of the Newton polygon of $\tilde{h}'$ whose slope is 
 $-(s+1)$. 

Consider the polynomial $\Theta(u)$ in one variable, defined by 
\[\Theta(\frac{z}{t^{s+1}}) = \frac{\init_{L'} 
\tilde{h}'}{t^{(s+1)(j_0-1)+k_0}}.\]

We argue as in Lemma \ref{lem:equid2} : the polynomial $\tilde{g}(u)$ 
has real roots at 0 and $c_{s+1}$, hence its derivative $\Theta(u)$ has at 
least one real root $d_{s+1} \in ]0,c_{s+1}[$. 

Hence, by (\ref{eq:init}), 
$g^{(k+1)}(\gamma_0)$ has at least one real root $v$ of the form 
$d_{s+1}t^{s+1} + h.o.t.$ with $0<d_{s+1} < c_{s+1}$.

This completes the proof of the Lemma. $\Box$
\medskip


\section{Reduction to the case when $\alpha$ and $\beta$ are curvettes with $\Gamma=\Z$ and $k_\alpha=k_\beta=R$}\label{Sec:Reduc} 

%
%

\begin{lem} \label{lem:generateurs}
 Let $\alpha,\beta \in \sper\ A$. Assume that $\alpha$ is not a specialization 
of $\beta$ and vice-versa. Then $<\alpha,\beta>$ is generated by all the 
elements $g \in A$ such that $g(\alpha)<0$ and $g(\beta)>0$.
\end{lem}

\noi Proof : Let $g_1,\ldots,g_s$ be a set of generators of $<\alpha,\beta>$ 
such that, for each $i$, we have $g_i(\alpha) \ge 0, g_i(\beta) \le 
0$. It suffices to show that for each 
$i \in \{1,\ldots,s\}$, there exist $h_{1 i}, h_{2i} \in A$ such that 
$h_{1i}(\alpha), h_{2i}(\alpha) >0$, $h_{1i}(\beta), h_{2i}(\beta) <0$ and $g_i 
\in (h_{1i},h_{2i})$. 

Since, by assumption, $\alpha$ is not a specialization of $\beta$, there exists 
$h_i \in A$ such that $h_i(\alpha) > 0$ and $h_i(\beta) \le 0$. Similarly, $\beta$ is not a specialization of $\alpha$, so there exists $k_i \in A$ such that $k_i(\alpha) \ge 0$, $k_i(\beta) <0$.  

Let $h_{1i}= g_i+h_i+k_i$ and $h_{2i}=2g_i+h_i+k_i$. So the desired $h_{1i}$ and 
$h_{2i}$ are constructed. $\Box$
\medskip

\noi \textbf{Notation} Let $A$ be a ring and $\alpha \in \sper\ A$. Let $P$ be a $\nu_\alpha$-ideal of $A$. We will denote by $P^+_\alpha$ the greatest $\nu_\alpha$-ideal of $A$, strictly contained in $P$.

\begin{prop}
 Let $A$ be a ring, $C$ an open set in $\sper\ A$ and 
 $\alpha,\beta \in C$ such that $\mathfrak{p} = \sqrt{<\alpha,\beta>}$ is a maximal ideal. There exist 
$\tilde{\alpha},\tilde{\beta} \in \sper\ A$ 
satisfying the following conditions :

(1) The value groups $\Gamma_{\tilde{\alpha}}$ of $\nu_{\tilde{\alpha}}$ and 
$\Gamma_{\tilde{\beta}}$ of $\nu_{\tilde{\beta}}$ are both isomorphic to $\Z$ and $k_\alpha \cong k_\beta \cong R$.

(2) $<\tilde{\alpha},\tilde{\beta}> \supset <\alpha,\beta>$, 
$<\tilde{\alpha},\tilde{\beta}>^+_{\tilde{\alpha}} \supset 
<\alpha,\beta>^+_\alpha$, 
$<\tilde{\alpha},\tilde{\beta}>^+_{\tilde{\beta}} \supset 
<\alpha,\beta>^+_\beta$.

(3) $\tilde{\alpha} \in C$ and $\tilde{\beta} \in C$.
\end{prop}

\noi Proof : Let $g_1,\ldots,g_s$ be as in the proof of Lemma \ref{lem:generateurs}. Let $g_1^\alpha, \ldots, g_{s^\alpha}^\alpha$ be a set of generators of $<\alpha,\beta>_\alpha^+$. Similarly let $g_1^\beta, \ldots, g_{s^\beta}^\beta$ be a set of generators of $<\alpha,\beta>_\beta^+$. Let $h_1,\ldots,h_r$ be a complete list of inequalities with appear in the definition of $C$.

Let $\pi:X \to \sper\ A$ be a sequence of blowings up with non singular centers having the following properties. In what follows, "prime" will denote strict transform under $\pi$. For example, $\alpha', \beta'$ are the strict transforms of $\alpha, \beta$ respectively. We require that :

1. the centers of $\alpha',\beta'$ in $X$ are disjoint;

2. for all indices $i$, the elements $g_i,g_i^\alpha,g_i^\beta,h_i$ define normal crossing subvarieties of $X$;

3. for all indices $i$, the sets $\{g'_i=0\}$, $\{g_i^{\alpha '}=0\}$, $\{g_i^{\beta '}=0\}$, $\{h_i=0$\} do not contain the center of $\alpha'$ or the center of $\beta'$.

Let $\tilde{\alpha}'$ be an $R$-semi-curvette with exponents in $\Z$ whose center is a sufficiently general closed point of the center of $\alpha'$ and similarly for $\tilde{\beta}'$. Then $\tilde{\alpha} := \pi(\tilde{\alpha'})$ and $\tilde{\beta} := \pi(\tilde{\beta'})$ satisfy the conclusions of the proposition. $\Box$

\section{Proof of the main theorem}\label{Sec:Proof}

We give a proof by contradiction. By Remark \ref{rek:centrecom}, we may assume that $\alpha, \beta$ have a common specialization $\xi$ and $f(\xi) = 0$. Assume that $\alpha$ and $\beta$ 
lie in 2 different connected components of $\{f\neq 0\}$.  Note that, because $f$ does not change sign between $\alpha$ and $\beta$, there are at least two real branches $f_1,f_2$ of $\{f=0\}$ over $D$, not necessarily distinct, containing $\xi$, between $\alpha$ and $\beta$. Hence there exists at least one real branch $h$
of $\{f'=0\}$ over $D$ between $\alpha$ and $\beta$. 
\medskip

Let $\theta = \min\{k >0 ; \ \nu_\alpha(f^{(k)}) < \nu_\alpha(f)\}$ (the set over which the minimum is taken is not empty because it contains $d$).
\medskip

Assume that there exists a real branch of 
$f^{(\theta)}$ over $D$ between $\alpha$ and $\beta$,
so $\alpha$ and $\beta$ lie in 
two different connected components of $C \cap \{f^{(\theta)} \neq 0\}$. Note that by the above this holds if $\theta=1$. By 
definition of $\theta$, $\nu_\alpha(f^{(\theta)}) < \nu_\alpha(f)$. 


If $f^{(\theta)}$ changes sign between $\alpha$ and $\beta$, 
$f^{(\theta)} \in <\alpha,\beta>$, which implies that $\nu_\alpha(f^{(\theta)}) \geq 
\mu_\alpha$ 
and $\nu_\beta(f^{(\theta)}) \geq \mu_\beta$. This contradicts the hypothesis 
$\nu_\alpha(f) \leq \mu_\alpha$ (and the same with $\beta$). Hence, 
$f^{(\theta)}$ does not change sign between $\alpha$ and $\beta$.

\begin{rek}
	If $\nu_\alpha(f^{(\theta)}) < \nu_\alpha(f)$, then $f^{(\theta)} \notin <\alpha,\beta>$ and 
	hence $\nu_\beta(f^{(\theta)}) < \mu_\beta$.
\end{rek}

Thus, if $\nu_\alpha(f^{(\theta)}) < \nu_\alpha(f)$, then the triple $(f^{(\theta)},\alpha,\beta)$ 
satisfies the hypothesis of the Theorem \ref{thm:main}. By induction on deg($f$), $\alpha$ and $\beta$ lie in the same connected component of $\{f^{(\theta)} \neq 0\} \cap C$. This is a contradiction. 
This completes the proof of the Theorem assuming there is a real branch of $f^{(\theta)}$ between $\alpha$ and $\beta$ (in particular in the case $\theta=1$).
\medskip

It remains to prove the existence of a real branch of $f^{(\theta)}$ between $\alpha$ and $\beta$. Since the Theorem has been proved in the cae $\theta=1$, we have reduced the problem to the case 
\begin{equation}
\nu_\alpha(f') \geq \nu_\alpha(f) \textbf{ and } \nu_\beta(f') \geq 
\nu_\beta(f) \label{eq:ineg2-2}
\end{equation}

We proceed by induction on $d=\deg(f)$. 
\smallskip

To prove the base of the induction, let us consider the cases $\deg(f)=1$ and 
$\deg(f)=2$. First, let $\deg(f)=1$. The set $\{f\neq 0\}$ has two connected 
components $C_1$ and $C_2$; up to 
interchanging $C_1$ and $C_2$, we have $f(C_1)>0$ and $f(C_2) <0$. 
Since $f$ does not change sign between $\alpha$ and $\beta$, $\alpha$ and 
$\beta$ lie in the same connected component of $\{f\neq 0\}$. 
\medskip

Next, let us consider the case $\deg(f)=2$. Assume that $\alpha$ and 
$\beta$ are in different components of $\{f\neq0\}$, aiming for 
contradiction. Write $f=(z-\phi_1)(z-\phi_2)$ where $f_1=z-\phi_1$, $f-2=z-\phi_2$ are branches. Since $(f,\alpha,\beta)$ are in good position, up to interchanging $f_1$ and $f_2$, we have $f_1(a) \le f_2(a)$ for all $a \in D$ and either equality holds for all $a \in D$ or strict inequality holds for all $a \in D$ (this is due to the choice of $D$, see Definition \ref{def:good} (2)). Hence $C \setminus \{f=0\} = C_{++} \coprod C_{--} \coprod C_{+-}$ where $C_{++} = \{f_1,f_2 >0\}$, $C_{--} =  \{f_1,f_2 <0\}$, $C_{+-} = \{f_1<0,f_2 >0\}$, where $C_{+-}$ may be empty; $f$ is strictly positive on $C_{++}$, $C_{--}$ and strictly negative on $C_{+-}$. Hence we may assume that $f_1(\alpha), 
f_2(\alpha) >0$ and $f_1(\beta), f_2(\beta) <0$. We have $f_1(\alpha) < f'(\alpha)$. Apply Lemma \ref{lem:conv2} to $f_1$ and $f'$. We obtain $\nu_\alpha(f') \leq \nu_\alpha(f_1)$. Since $\nu_\alpha(f_2) >0$, we have $\nu_\alpha(f') < \nu_\alpha(f)$. Which is a contradiction. 
\medskip

The base of the induction is proved, let us now prove the induction step. 
\medskip

Let $f=\prod_{j=1}^d f_j$, let $f'=\prod_{j=1}^{d-1}h_j$ where $f_j$ are $\alpha_0$-branches of $f$ and $h_j$ are $\alpha_0$-branches of $f'$.  
Let $h_\alpha$ be an $\alpha$-privileged branch of $f'$. Then, since $\nu_\alpha(f') 
\ge \nu_\alpha(f)$, we have  
\begin{equation}
 \nu_\alpha(Im\, h_\alpha) \ge \nu_\alpha(h_\alpha) > \nu_\alpha(f_j), \ \text{ 
for } 1 \le j \le d.
\end{equation}
The first inequality is by Lemma \ref{lem:reim} and the second is just (\ref{eq:inegval}).
\medskip

\noi \textbf{Claim :} Let $i \in \{0,\ldots, [\frac{\theta}{2}] \}$. There exist  :

- branches over $D$, $g_{2i,1}, g_{2i,2}, $ of $f^{(2i)}$ and $\tilde{h}_{2i+1}$ of $f^{(2i+1)}$,

-  a branch $h_{2i+1,\alpha}$of $f^{(2i+1)}$ over a suitable neighbourhood $U_\alpha$ of $\alpha_0$ and a branch $h_{2i+1,\beta}$ of $f^{(2i+1)}$, over a suitable neighbourhood $U_\beta$ of $\beta_0$ 

such that : 

(1) $g_{2i,1}, g_{2i,2},\tilde{h}_{2i+1}$ are real and separate $\alpha$ 
from $\beta$ and \begin{align} \label{eq:separate}
        0 < g_{2i,1}(\alpha) < \tilde{h}_{2i+1}(\alpha) < g_{2i,2}(\alpha) 
&\text{ or } 0 > g_{2i,1}(\alpha) > \tilde{h}_{2i+1}(\alpha) > 
g_{2i,2}(\alpha)  \\  g_{2i,1}(\beta) > \tilde{h}_{2i+1}(\beta) > 
g_{2i,2}(\beta) > 0 &\text{ or } g_{2i,1}(\beta) < \tilde{h}_{2i+1}(\beta) < 
g_{2i,2}(\beta) < 0 
                 \end{align}

(2) for $i>0$ we have \begin{align}
                       \nu_\alpha(g_{2i,2}(\alpha)) \le 
\nu_\alpha(g_{2i,1}(\alpha)) &=\nu_\alpha(\tilde{h}_{2i-1}(\alpha)) < 
\nu_\alpha(h_{2i-1,\alpha}(\alpha)) \\
\nu_\beta(g_{2i,1}(\beta)) \le \nu_\beta(g_{2i,2}(\beta)) 
&=\nu_\beta(\tilde{h}_{2i-1}(\beta)) < \nu_\beta(h_{2i-1,\beta}(\beta));
                      \end{align}

(3) $h_{2i+1,\alpha}$ is an $\alpha$-privileged branch of $f^{(2i+1)}$ (and 
similarly for $h_{2i+1,\beta}$); 

(4) \begin{align}
    \nu_\alpha(g_{2i,1}(\alpha)) &\geq \nu_\alpha(\tilde{h}_{2i+1}(\alpha)) = 
\nu_\alpha(g_{2i,2}(\alpha)) \\ \nu_\beta(g_{2i,2}(\beta)) &\geq 
 \nu_\beta(\tilde{h}_{2i+1}(\beta)) = \nu_\beta(g_{2i,1}(\beta)).
    \end{align}
%
%
\medskip

\noi Proof of Claim : We use induction on $i$. First let 
$i=0$. We put 
$g_{0,1}=f_1, g_{0,2}=f_2,\tilde{h}_1=h$ with the convention that $f_1$ is 
between $\alpha$ and $f_2$ and $f_2$ is between $\beta$ and $f_1$.
Fix an 
$\alpha$-privileged branch $h_{1,\alpha}$ of $f'$ and 
similarly for $\beta$. 
Statements (1), (3) are clear, (4) follows from Lemmas \ref{lem:conv2} and 
\ref{lem:equid2}. (2) is vacuously true. 
\smallskip

Let $0<i\leq [\frac{\theta}{2}]$. We assume the Claim holds up to 
$i-1$. 

As $i \le 
[\frac{\theta}{2}]$, we have $\nu_\alpha(f^{(2i-1)}) \ge \nu_\alpha(f)$. 
Using (\ref{eq:inegval}), we have $\nu_\alpha(h_{2i-1,\alpha}(\alpha)) > 
\nu_\alpha(g_{0,1}(\alpha))$. By (\ref{eq:separate}) at step $(i-1)$, we may apply 
Lemma \ref{lem:conv2} to 
$g_{0,1},g_{2i-2,2}$. We obtain $\nu_\alpha(g_{0,1}(\alpha)) \ge 
\nu_\alpha(g_{2i-2,2}(\alpha))$. According to (4) at step $i-1$, we have $\nu_\alpha(g_{2i-2,2}(\alpha)) = \nu_\alpha(\tilde{h}_{2i-1}(\alpha))$. So 
\begin{equation} \label{eq:compar}
 \nu_\alpha(\tilde{h}_{2i-1}(\alpha)) < \nu_\alpha(h_{2i-1,\alpha}(\alpha)). 
\end{equation}

Let $g_{2i,1}(\alpha)$ be a real factor of $f^{(2i)}(\alpha)$ between $\tilde{h}_{2i-1}(\alpha)$ and $h_{2i-1,\alpha}$ whose existence is given by Lemma \ref{lem:equid3}. Let $g_{2i,1}$ be a real branch of $f^{(2i)}$ associated to $g_{2i,1}(\alpha)$ by  Remark \ref{rek:bijection}. 
Similarly, let $g_{2i,2}$ be a real 
branch of $f^{(2i)}$ over $D$ between $\tilde{h}_{2i-1}$ 
and $h_{2i-1,\beta}$. 

Let $\tilde{h}_{2i+1}$ be a real branch of $f^{(2i+1)}$ over $D$ between $g_{2i,1}$ and 
$g_{2i,2}$ whose existence is guaranteed by Rolle's theorem. 

Let $h_{2i+1,\alpha}$ be an $\alpha$-privileged branch of $f^{(2i+1)}$ and 
similarly for $\beta$. Statements (1) and (3) are clear. (2) follows from 
(\ref{eq:compar}) and Lemma \ref{lem:conv2}. 
Now (4) follows from (1) and Lemma \ref{lem:equid2}. 

This completes the proof of 
the Claim. $\Box$
\medskip

We continue with the proof of Theorem \ref{thm:main}. Let 
$i=[\frac{\theta}{2}]$. 
Apply the Claim to this $i$. By the Claim, there exists a real branch of 
$f^{(\theta)}$ over $D$ between $\alpha$ and $\beta$.
This completes the proof of the Theorem.
$\Box$ \medskip

	Given a monic polynomial $g \in B[z]$, write $g = \prod_{j=1}^s g_j^{a_j}$, where each $g_j$ is a linear polynomial over the algebraic closure of the field of fractions of $B$. 
	
	\begin{deft}
		A partial reduction  of $g$ is a polynomial of the form $\tilde{g} = \prod_{j=1}^{s} g_j^{a'_j}$, where $0< a'_j \le a_j$ for all $j$. A derivative reduction of $g$ is a polynomial of the form $\tilde{g}'$ where $\tilde{g}$ is as above.
	\end{deft}
	
	\begin{deft}
		$(f,\alpha,\beta)$ are in generalized good position if there exists a sequence of polynomial having the following properties :
		
		(1) For each $i \in \{0,\ldots,k\}$, the polynomial $f^{(i+1,*)}$ is a derivative reduction of $f^{(i*)}$;
		
		(2) $\nu_\alpha(f) \le \nu_\alpha(f^{(1*)}), \ldots,\nu_\alpha(f^{(i-1,*)})$ and $\nu_\alpha(f) > \nu_\alpha(f^{(i*)})$;
		
		(3) The number of real roots of each of $f,f^{(1*)}, \ldots,f^{(i*)}$ counted with or without multiplicity is constant over $D$.
	\end{deft}

\begin{rek} \label{rek:extension}
	With the same proof, the main theorem holds for $(f,\alpha,\beta)$ in generalized good position.
\end{rek}

\end{document}